\def\al{abstract linkage}

\def\rl{cabled linkage}
\def\ul{{\cal L}}
\def\fdf#1{{\it #1}}
\def\conf{{\cal C}}
\def\vert{{\cal V}}
\def\edge{{\cal E}}
\def\qas{quasialgebraic set}
\def\KM{Kapovitch and Millson}
\def\mand{\ {\rm and}\ }
\def\mif{\ {\rm if}\ }
\def\qf{quasifunctional}
\def\rd{restricted domain}
\def\cite#1{{\bf[#1]}}

\def\compose{\circ}
\def\PC{Peaucellier}
\def\rh#1#2{\rho _{\ul #1,#2}}
\def\sec#1{\goodbreak\bigbreak \centerline{\bf #1} \medskip}
\def\capshun#1#2{\centerline{{\bf Figure #1:} #2}}
\def\mapright#1{\smash{\mathop{\longrightarrow}\limits^{#1}}}
\def\mapdown#1{\Big\downarrow\rlap{$\vcenter{\hbox{$\scriptstyle#1$}}$}}
\def\commentout#1{}
\def\proof{\par{\bf\noindent Proof: }}
\def\qed{{\hskip 0pt plus 1filll}
 \vbox{\hrule height 4pt width 4pt}\hfil \par \medbreak}
 %  change these depending on whether function vertices can be duplicated
\def\ifdupv#1{#1}
\def\ifnotdupv#1{}

%\commentout{
\catcode`\@=11
\font\tenmsb=msbm10
\font\sevenmsb=msbm7
\font\fivemsb=msbm5
\newfam\msbfam
\textfont\msbfam=\tenmsb  \scriptfont\msbfam=\sevenmsb
  \scriptscriptfont\msbfam=\fivemsb
\def\Bbb{\ifmmode\let\next\Bbb@\else
 \def\next{\errmessage{Use \string\Bbb\space only in math mode}}\fi\next}
\def\Bbb@#1{{\Bbb@@{#1}}}
\def\Bbb@@#1{\fam\msbfam#1}
\catcode`\@=12
%}
%\def\Bbb#1{{\bf #1}}

\magnification\magstep1

\centerline{\bf Planar Linkages and Algebraic Sets}
\medskip
\centerline{\bf Henry C. King}
\bigskip

An \fdf{\al} is a finite graph $L$ with a positive number
$\ell (\overline {vw})$ assigned to each edge $\overline {vw}$.
A \fdf{planar realization} of an \al\ $(L,\ell )$ is a mapping
$\varphi $ from the vertices of $L$ to ${\Bbb C}$ so that 
$|\varphi (v)-\varphi (w)|=\ell (\overline{vw})$ for all  edges $\overline {vw}$.
We will investigate the topology of the space of planar realizations
of linkages.

You may think of an \al\ as an ideal mechanical device
consisting of a bunch of stiff rods (the edges) with length
given by $\ell $ and sometimes attached at their ends
by rotating joints. A planar realization is some way of
placing this linkage in the plane.

If $L$ is a finite graph, we let $\vert(L)$ denote the set of vertices of $L$
and let $\edge(L)$ denote the set of edges of $L$.

We will often wish to fix some of the vertices of a linkage whenever we
take a planar realization.
So we say that a \fdf{linkage} $\ul$ is a foursome $(L,\ell ,V,\mu )$ where
 $(L,\ell )$ is an \al, $V\subset \vert(L)$ is a subset of its vertices, and $\mu \colon V\to {\Bbb C}$.
 So $V$ is the set of fixed vertices and $\mu $ tells where to fix them.
The configuration space of realizations is defined by:
$$\conf(\ul)=\{\,\varphi \colon \vert(L)\to {\Bbb C}\,\mid \, \varphi (v)=\mu (v)\mif\ v\in V, \ {\rm and}\
   |\varphi (v)-\varphi (w)|=\ell (\overline{vw})\ {\rm for\ all\ edges}\ \overline {vw}\,\}$$
If $w$ is a vertex of $\ul$, we have a map $\rh{}w\colon \conf(\ul)\to {\Bbb C}$
 given by $\rh{}w(\varphi )=\varphi (w)$, the position of the vertex $w$.
In this paper we will look at configuration spaces $\ul$
and also the image of $\rh{}w$.
   
The literature on linkages goes back over 300 years.
Moreover, humans have been using linkages for thousands of years
in various mechanical devices.

I recall first hearing of planar linkages in a talk by Thurston
at IAS, probably in the mid 1970s.
I recall three results from this talk:
\item{1)} Any algebraic set is isomorphic to an algebraic set
given by quadratic polynomials.
\item{2)} You can construct a linkage which will ``sign your name",
i.e., there is a linkage $\ul$ and vertices $v_1,\ldots ,v_k$ so that
$\cup _{i=1}^k \rh{}{v_i}(\conf(\ul))$ is an arbitrarily close approximation
of your signature.
\item{3)} Some result about realizing any compact smooth manifold
as a configuration space of a linkage, which result I do not recall at all precisely.

As far as I can tell, Thurston never wrote these results up, so 3)
must remain vague.
Occasionally since then I have been contacted by an engineer interested
in these results, but I could not recall anything about Thurston's
proof so I could not help them.
Then recently John Millson started asking me lots of questions on
real algebraic sets.  He and M.~Kapovich were writing up proofs
of the results 2) and 3) above.
In the course of doing so, they discovered and solved
some problems overlooked by previous literature.
In this paper I will give a proof of 2) and 3) based upon
the Kapovich-Millson proof in \cite{KM}.  In particular
I will give proofs of the following:

\proclaim{Theorem 1}. (Thurston?,Kapovich-Millson)
If $X$ is a compact real algebraic set,
then there is a linkage $\ul$ so that $\conf(\ul)$
is the union of a number of copies of $X$, 
each analytically isomorphic to $X$.
In fact there is a polynomial map from $\conf(\ul)$
to $X$ giving an analytically trivial covering of $X$.

\proclaim{Theorem 2}. (Thurston,Kapovich-Millson)
If $\alpha \colon [a,b]\to {\Bbb C}$ is any polynomial map, there is a linkage $\ul$
and a vertex $v$ of $\ul$ so that $\rh{}v$ traces out the curve
$\alpha $, i.e., $\alpha ([a,b])=\rh{}v(\conf(\ul))$.

%\proclaim{Theorem 2.1}.(King?)
%If $\alpha _i\colon [a_i,b_i]\to {\Bbb C}$ are polynomial maps, $i=1,\ldots ,k$, 
%there is a linkage $\ul$
%and a vertex $v$ of $\ul$ so that $\rh{}v$ traces out the union of the curves
%$\alpha _i$, i.e., $\Bigcup _{i=1}^k \alpha _i([a_i,b_i])=\rh{}v(\conf(\ul))$.

\proclaim{Corollary 2.1}. (Thurston)
There is a linkage which signs your name.

The corollary follows from Theorem 2 because you can 
break up your signature into $k$ segments, each of which is
the image of a smooth curve, and then approximate each
smooth curve by a polynomial curve.
If you wish, you may even do this approximation in a way which
preserves any cusps in your signature.
By being careful, you can even show there is a linkage $\ul$ which
signs your name using just one output vertex $v$.
So $\rh{}{v}(\conf(\ul))$ is a close approximation of your signature,
no matter how many strokes are needed, see \cite{K}.

I should define some terms here.
A \fdf{real algebraic set} is the set of solutions of a collection of
real polynomial equations in ${\Bbb R}^n$.
Note that if we view ${\Bbb C}$ as ${\Bbb R}^2$, any configuration
space $\conf(\ul)$ is an algebraic set in ${\Bbb R}^{2\vert(\ul)}$,
since it is the solutions of the polynomial equations
$|y_i-y_j|^2=\ell (ij)^2$ and $y_i=z_i$ for $i\in V$.
Moreover these equations are quadratic at worst
which is the reason for Thurston's result 1) above
which says that quadratic algebraic sets are
not at all special as algebraic sets go.
\footnote{${}^1$}{I will illustrate the proof of 1) with an example.
If we make a new variable $u$ and set $u=x^2$,
then the cubic $x^3$ becomes the quadratic $ux$.}

We will use two notions of isomorphism.
If $X\subset {\Bbb R}^n$ and $Y\subset {\Bbb R}^m$ then we say a 
homeomorphism $f\colon X\to Y$ is an \fdf{isomorphism}
if $f$ and $f^{-1}$ are both restrictions of rational maps,
for example polynomials.
We say $f$ is an \fdf{analytic isomorphism}
if $f$ and $f^{-1}$ are both restrictions of analytic maps,
i.e., maps locally given by power series.
So any isomorphism is analytic, but the converse is not true.

\sec{Cabled Linkages}

Note that Theorem 1 is somewhat unsatisfying,
since it does not give a complete characterization of all
compact configuration spaces, but only a characterization
up to finite analytically trivial covers.

So we introduce the following generalization of linkage
which will allow a complete characterization of
configuration spaces.

A \fdf{\rl}\ is a quintuple $(L,\ell ,V,\mu ,F)$ where
$(L,\ell ,V,\mu )$ is a linkage and $F\subset \edge(L)$.
We will think of the edges in $F$ as being flexible
rather than rigid.
A physical model for such a \rl\ would
 have the edges in $\edge(L)-F$ be rigid rods as before
 but the edges in $F$ are just ropes or cables
connecting two vertices.
Thus in a planar realization, two vertices connected by an edge $e$ in
$F$ would only be constrained to have distance $\le \ell (e)$.
If $F$ is empty we get a classical linkage.
The configuration space is given by:
\commentout{
$$\eqalign{
\conf(\ul)&=\{\,\varphi \colon \vert(L)\to {\Bbb C}\,\mid \, \varphi (v)=\mu (v)\ {\rm if}\ v\in V,\cr
 &\quad \quad\quad \quad\quad \quad |\varphi (v)-\varphi (w)|\le \ell (\overline{vw})\mif\ \overline {vw}\in F \mand\ \cr
   &\quad\quad\quad \quad\quad \quad |\varphi (v)-\varphi (w)|=\ell (\overline{vw})\mif\ \overline {vw}\in \edge(L)-F\,\}
}$$
}
$$\conf(\ul)=\left\{ \varphi \colon \vert(L)\to {\Bbb C} \Biggm|
\matrix{
&\varphi (v)=\mu (v)\ {\rm if}\ v\in V,\cr
  &|\varphi (v)-\varphi (w)|\le \ell (\overline{vw})\mif\ \overline {vw}\in F \mand\ \cr
   &|\varphi (v)-\varphi (w)|=\ell (\overline{vw})\mif\ \overline {vw}\in \edge(L)-F\cr
}
\right\}$$

From now on, the word linkage will refer to a \rl.
If we wish to refer to a linkage without any flexible edges,
we will call it a classical linkage.
   
 Because of the inequalities, the configuration space 
of a \rl\ may no longer be
 a real algebraic set.
 However it is something which I will call a \qas.
 I define a \fdf{\qas}\ to be a subset of ${\Bbb R}^n$ of the form
 $$\{\,x\in {\Bbb R}^n\,\mid\, p_i(x)=0,\, i=1,\ldots ,k \mand\ q_j(x)\ge 0,\, j=1,\ldots ,m\,\}$$
 for some polynomials $p_i$ and $q_j$.
 Using \rl s, we get the following complete characterization
 of configuration spaces of \rl s.
 
\proclaim{Theorem 3}. (King)
If $X$ is a compact \qas,
then there is a \rl\ $\ul$ so that $\conf(\ul)$
is analytically isomorphic to $X$.
In fact there is a polynomial map from $\conf(\ul)$
to $X$ giving this isomorphism.

\proclaim{Corollary 3.1}. (King)  
Up to analytic isomorphism, the set of configuration spaces
of \rl s is exactly the set of spaces of the form
$X\times {\Bbb R}^{2k}$ where $X$ is a compact \qas.

So, for example using \cite{AK} and \cite{AT}
any compact smooth manifold with boundary is diffeomorphic to 
some $\conf(\ul)$ and any compact PL manifold with boundary 
is homeomorphic to  some $\conf(\ul)$.

\sec{Functoriality of $\conf(\ul)$}

Let $\ul'\subset \ul$ be a sublinkage.
This means that $L'\subset L$, $\ell '=\ell |_{\edge(L')}$, 
$V'\subset V$, $\mu '=\mu |_{V'}$,
and $F'=F\cap \edge(L')$.
Then we have a natural map $\rh{}{\ul'}\colon \conf(\ul)\to \conf(\ul')$
obtained by restriction, i.e., $\rh{}{\ul'}(\varphi )=\varphi |_{\vert(L')}$.
Note that if $L'$ is a single point $v$ and $V'$ is empty, then
$\conf(\ul')={\Bbb C}$ and with this identification we have $\rh{}{\ul'}=\rh{}v$.

If we wished, we could generalize this by forming a category of linkages
and linkage maps and noting that $\conf$ is a contravariant functor.
But we really only have a need to look at inclusion maps,
\ifdupv{or sometimes simple quotients,}
so we refrain from this generality.

\proclaim{Lemma 3.2}.
If $\ul'\subset \ul$ and $\ul''\subset \ul$ are two sublinkages then we have
a natural identification of $\conf(\ul'\cup \ul'')$ with the fiber product
of the restriction maps $\rh{'}{\ul'\cap \ul''}\colon \conf(\ul')\to \conf(\ul'\cap \ul'')$
and $\rh{''}{\ul'\cap \ul''}\colon \conf(\ul'')\to \conf(\ul'\cap \ul'')$.

% put in commut diagram
$$\matrix{ 
\conf(\ul'\cup \ul'') & \mapright{} & \conf(\ul'') \cr
\mapdown{}           &             &   \mapdown {\rh{''}{\ul'\cap \ul''}} \cr
\conf(\ul')     & \mapright{\rh{'}{\ul'\cap \ul''}}  & \conf(\ul'\cap \ul'') \cr
}$$

\proof
This is because a planar realization of $\ul'\cup \ul''$ is just
a planar realization of $\ul'$ and a planar realization of $\ul''$
which happen to agree on $\ul'\cap \ul''$.
Thus 
$$\conf(\ul'\cup \ul'')=\{(\varphi ',\varphi '')\in \conf(\ul')\times \conf(\ul'') \mid 
\rh{'}{\ul'\cap \ul''}(\varphi  ')=\rh{'}{\ul'\cap \ul''}( \varphi '')\}$$
is the fiber product.
\qed

As a consequence of Lemma 3.2, the configuration space of
the disjoint union of two linkages is the product of their
configuration spaces.

\proclaim{Lemma 3.3}.
Suppose  $\ul$ is a \rl\ with no fixed vertices.
Form a new $\ul'$ from $\ul$ by fixing exactly
one vertex in each connected component of $\ul$.
If $\ul$ has $k$ connected components, then
$\conf(\ul)$ is isomorphic to $\conf(\ul')\times {\Bbb C}^k$.

%Suppose $\ul=(L,\ell ,V,\mu ,F)$ is a \rl\ with $V$ empty.
%Pick any $V'\subset \vert(L)$ so that $V'$ has exactly one vertex 
%in each connected component of $L$.
%Choose any map $\mu '\colon V'\to {\Bbb C}$.
%Consider the new linkage $\ul'=(L,\ell ,V',\mu ',F)$.
%Note then that if $L$ has $k$ components then
%$$\conf(\ul)=\conf(\ul')\times {\Bbb C}^k$$
\proof
By the above remark on the configuration space of a 
disjoint union, it suffices to show this when $\ul$ is connected,
 so $k=1$.
For the map from right to left,
take any $\varphi '\in \conf(\ul')$ and any $z\in {\Bbb C}$.
We then get a $\varphi \in \conf(\ul)$ by letting  $\varphi $ 
be the composition of $\varphi '$ and translation by $z$.
\qed

So if we want to understand the topology configuration spaces,
it suffices to consider only those \rl s with at least one
fixed vertex in each connected component of $L$.
Note that such configuration spaces must be compact
since the distance between any two vertices in the
same component of $L$ is bounded by the sum of the
lengths of the edges of a path connecting them.

\proclaim{Lemma 3.4}. Let $\ul$ be a \rl\ and let
$v_1,\ldots ,v_m$ be vertices of $\ul$ which are not fixed.
Let $\ul'$ be obtained from $\ul$ by fixing the vertices $v_1,\ldots ,v_m$
to be at the points $z_1,\ldots ,z_m$.
Let $p\colon \conf(\ul)\to {\Bbb C}^m$ be the map $(\rh{}{v_1},\ldots ,\rh{}{v_m})$.
Then $\conf(\ul')=p^{-1}(z_1,\ldots ,z_m)$.

%To prove this, note that $\ul'=\ul\cup \ul''$ where $\ul''$
%is the linkage with no edges and vertices $v_1,\ldots ,v_m$,
%all of which are fixed (to the points  $z_1,\ldots ,z_m$).
%Note that $\conf(\ul'')$ is a point.
%Also $\ul\cap \ul''$ is the linkage with no edges and with
%vertices $v_1,\ldots ,v_m$,
%none of which are fixed.
%So $\conf(\ul\cap \ul'')={\Bbb C}^m$ and the image of
%$\rh{}{\ul\cap \ul''}\colon \conf(\ul'')\to \conf(\ul\cap \ul'')$ is the point $(z_1,\ldots ,z_m)$.
%Hence we get the result from the fiber product.
%
%Alternatively, just do this directly.
\proof
Note that $\conf(\ul')\subset \conf(\ul)$ and it must be exactly those
$\varphi $ with $\rh{}{v_i}(\varphi )=\varphi (v_i)=z_i$.
The lemma follows.
\qed

\proclaim{Lemma 3.5}. Let $\ul$ be a \rl\ and let
$v_1,\ldots ,v_m$ be vertices of $\ul$ which are not fixed.
Let $\ul'$ be obtained from $\ul$ by 
adding new vertices $u_1,\ldots ,u_m$ and new flexible edges
$u_iv_i$ of length $b_i$.
Fix the vertices $u_i$ to points $z_i\in {\Bbb C}$.
Let $p\colon \conf(\ul)\to {\Bbb C}^m$ be the map $(\rh{}{v_1},\ldots ,\rh{}{v_m})$.
Then $\conf(\ul')$
is isomorphic to $p^{-1}(\{w\in {\Bbb C}^m \mid b_i\ge |w_i-z_i|, i=1,\ldots ,m\})$.

\proof
The inclusion $\ul\subset \ul'$ gives a map $\alpha \colon \conf(\ul')\to \conf(\ul)$.
%Note that $\alpha $ is an imbedding since the vertices of $\ul'$ which are
%not in $\ul$ are all fixed, so $\alpha $ has a trivial inverse.
%Also, $\alpha (\conf(\ul'))$ must be exactly those
%$\varphi $ with $|\rh{}{v_i}(\varphi )-z_i|\le b_i$.
Let $Y=p^{-1}(\{w\in {\Bbb C}^m \mid b_i\ge |w_i-z_i|, i=1,\ldots ,m\})$.
We have a map $\beta \colon Y\to \conf(\ul')$ defined by
$\beta (\varphi )(v)=\varphi (v)$ for $v$ a vertex of $\ul$ and
$\beta (\varphi )(u_i)=z_i$.
Note that $\alpha (\conf(\ul'))\subset Y$ and $\beta $ is the inverse of
$\alpha \colon \conf(\ul')\to Y$.
\qed

The following is immediate from the definitions.

\proclaim{Lemma 3.6}.
If $\ul'\subset \ul$ is a  sublinkage then
The map $\rh{}{\ul'}\colon \conf(\ul)\to \conf(\ul')$ is an (analytic)
isomorphism if and only if it is onto and the position
$\varphi (v)$ of each vertex $v$ of $\ul$ is a rational (resp.~analytic) function
of the positions $\varphi (w_i)$ of the vertices $w_i$ in $\ul'$.
More generally, if $Z\subset \conf(\ul)$ then the restriction
$\rh{}{\ul'}|_Z\colon Z\to \rh{}{\ul'}(Z)$ is an (analytic) isomorphism if and only if 
for $\varphi \in Z$, the position
$\varphi (v)$ of each vertex $v$ of $\ul$ is a rational (resp.~analytic) function
of the positions $\varphi (w_i)$ of the vertices $w_i$ in $\ul'$.

\proclaim{Lemma 3.7}.
\ifdupv{
Let $\ul$ be a linkage and suppose $v$ and $w$ are two vertices of $\ul$.
Suppose that whenever there are edges $vu$ and $wu$ to the same vertex $u$,
 that $\ell (vu)=\ell (vw)$.
Suppose also that there is no edge $vw$.
Then we may form a linkage $\ul'$ from $\ul$ 
by identifying the vertices $v$ and $w$,
and identifying any edges $vu$ and $wu$.
Moreover there is a natural identification of $\conf(\ul')$
with $\{\varphi \in \conf(\ul)\mid \varphi (v)=\varphi (w) \}$.
}

\sec{How to put a rotating joint in the middle of an edge}

In our \rl s, rotating joints between two edges only occur
at the ends of an edge.
In real life, a mechanical linkage may have a joint in the middle.
We may simulate such a joint as follows:

% insert picture fig 1
\midinsert\vskip 1.1in
\includegraphics{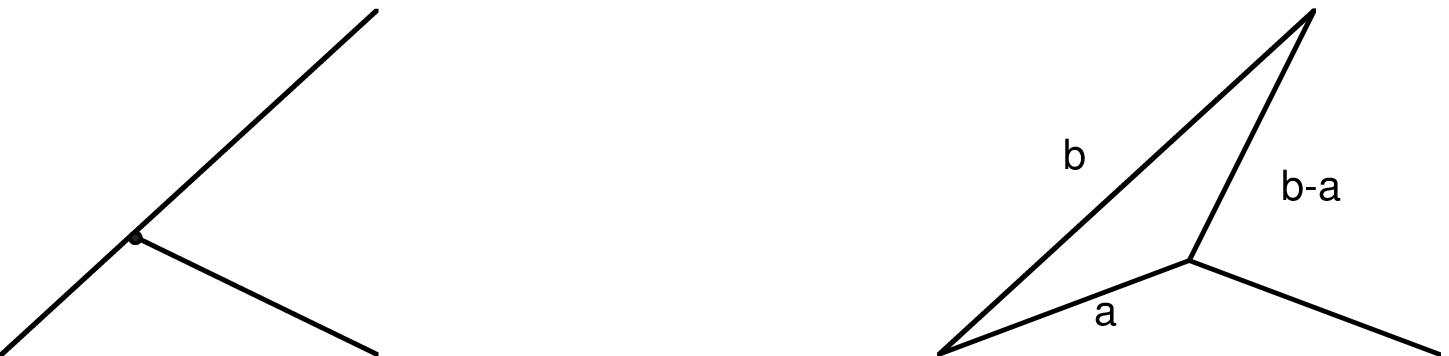}
\capshun 1 { How to put a joint in the middle of an edge}
\endinsert

Thus when drawing linkages, it is allowable to draw
a joint in the middle of an edge.

Note, this paragraph is optional.
If you actually constructed the linkage in Figure 1, you would
find it somewhat spongy, the position of the middle point
is quite sensitive to small errors in length.
This is reflected algebraically in the following.
Suppose you
look at the configuration space as a scheme, i.e., you
focus on the polynomial equations defining $\conf(\ul)$
rather than the point set of their solutions.
Then in the above linkage, the configuration scheme
is not reduced.
(This is because of the nontransversality of the
equations specifying the middle point.)
In \cite{KM}, where they take a scheme-theoretic point of view,
they get around this by essentially allowing joints in the
middle of an edge and modifying the equations accordingly.
One could also add a stiffening truss as shown below
in Figure 2
and the resulting configuration scheme is reduced.
In Figure 2,  $c^2=a^2+d^2$ and $e^2=d^2+(b-a)^2$.

% insert Fig 2
\midinsert\vskip 1.2in
\includegraphics{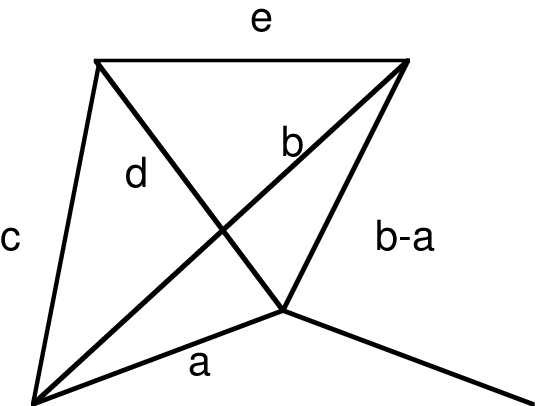}
\capshun 2 {Stiffening a middle joint }
\endinsert

\sec{The Square Linkage (actually a rhombus)}

Let us now look at a very simple linkage and describe its
configuration space.
Consider the square linkage in Figure 3 below.
As we will see soon, this is a sublinkage of a number
of useful linkages.
Let us try to find its configuration space if we
 fix the two bottom vertices $A$ and $B$.
It is tempting to believe this configuration space is a circle,
that once the position of $C$ is decided then the position
of vertex $D$ is determined.
However, Kapovitch and Millson noticed that in fact there
are other degenerate realizations, for example D could map
to the same point as A, leaving $C$ to rotate.
Or $C$ could map to the same point as $B$ and $D$
could rotate.
So in fact the configuration space is the union of three circles
each pair intersecting in a single point.
A similar problem occurs for a rectangle which is not a square.
Its configuration space is the union of two circles which 
intersect in two points.

% insert fig 3
\midinsert\vskip 1.5in
\includegraphics{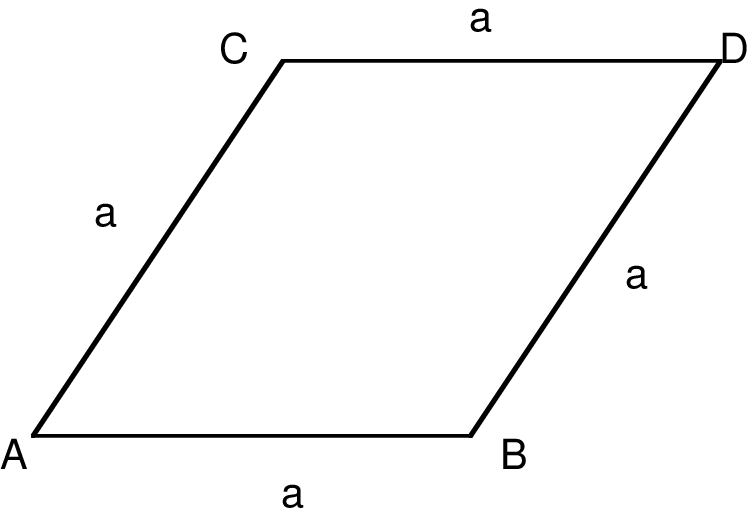}
\capshun 3 { The Square Linkage, with an unexpected Configuration Space}
\endinsert

In the linkages we wish to use, we will want to get rid of all
these degenerate realizations.
\KM\ did this by rigidifying the square.
They added another edge between the midpoints of two opposite
sides.
Note that for this new linkage, the degenerate realizations will not occur.
In order not to clutter up our pictures, we will subsequently draw this
extra edge as a dotted line.

% Figure 4 here
\midinsert\vskip 1.5in
\includegraphics{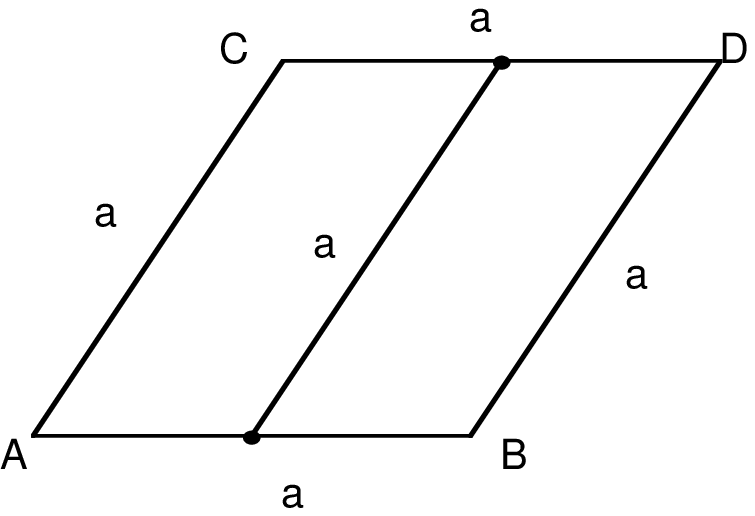}
\capshun 4 { The Rigidified Square Linkage}
\endinsert

\sec{Functional Linkages}

Let us now look at a very useful linkage, the pantograph.
The pantograph has been used as a mechanical device for centuries
 (without the rigidifying extra edge).
It can do a number of useful things.
For example, if $A$ is fixed at 0 and $B$ is at some point $z$,
then $C$ will be at the point $(1+c)z$.
This is an example of a functional linkage, since it 
can be used to evaluate the function $z\mapsto \lambda z$.

% Figure 5 pantograph here
\midinsert\vskip 1.5in
\includegraphics{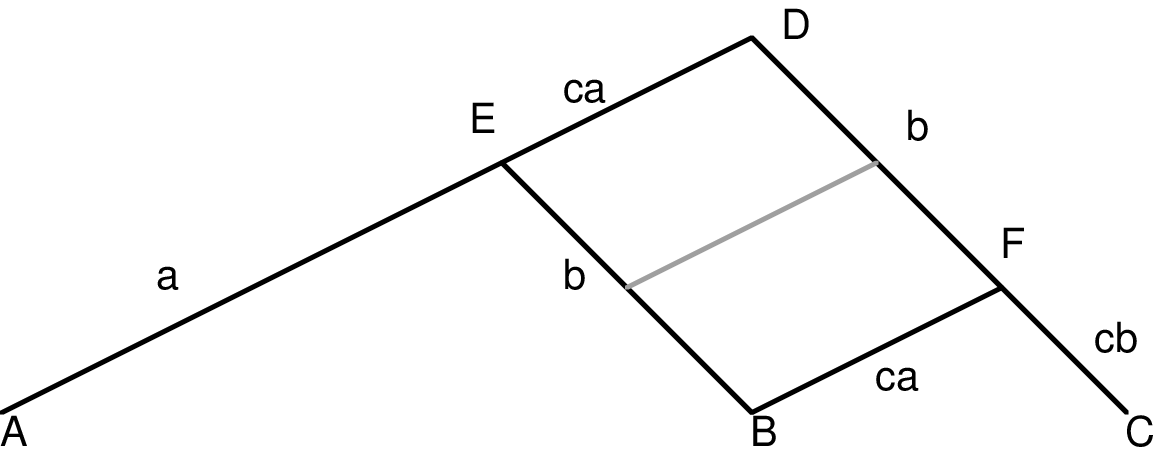}
\capshun 5 { The Pantograph}
\endinsert

Notice however that $B$ cannot move to every point of the plane.
Since $A$ is fixed at 0,
the point $B$ is constrained to lie in the annulus between the
circles of radius $a-b$ and $a+b$.
Moreover, for any position of the point $B$ in the interior of this annulus,
 there will be exactly two configurations depending on which
 side of the line $AB$ the point $D$ lies.
 This generalizes to the following:

A linkage $\ul$ is \fdf{\qf\ } for a map $f\colon {\Bbb C}^n\to {\Bbb C}^m$ 
 if there are 
 vertices $w_1,\ldots ,w_n$ and $v_1,\ldots ,v_m$ of $\ul$ so that
 if $p\colon \conf(\ul)\to {\Bbb C}^m$ is $(\rh{}{v_1},\ldots ,\rh{}{v_m})$ and
 $q\colon \conf(\ul)\to {\Bbb C}^n$ is $(\rh{}{w_1},\ldots ,\rh{}{w_n})$ then $p=f\compose q$.
The set $q(\conf(\ul))$ is called the \fdf{domain} of the \qf\ linkage.
We call $p$  the input map and call $q$ the output map.

If in addition, there is a $U\subset q(\conf(\ul))$ 
so that the restriction $q|\colon q^{-1}(U)\to U$ is an analytically trivial
 covering map, we say that 
$\ul$ is \fdf{functional} for $f$ 
 with \rd\ $U$.

%Moreover, if the restriction $q|\colon q^{-1}(U)\to U$ is an analytic isomorphism, we say that $\ul$ is \fdf{strongly functional} for a $f$ 
Moreover, if $q\colon \conf(\ul)\to q(\conf(\ul))$ is an analytic isomorphism
and $U\subset q(\conf(\ul))$ we say that $\ul$ is \fdf{strongly functional} for $f$ 
  with \rd\ $U$.

 We call $w_1,\ldots ,w_n$ the input vertices and call
 $v_1,\ldots ,v_m$ the output vertices.
 \ifdupv{Repetitions of vertices are allowed,
 although they are not necessary for the results in this paper.}
\ifnotdupv{
 We require all input vertices to be distinct,
 and all output vertices to be distinct.
 We will, however, allow an input vertex to also be an output vertex.
 % If we did not require distinction, we would have trouble
 % composing two functional linkages
 % composition would entail identifying vertices
 % but in a linkage you cannot just arbitrarily identify
}

So if $\ul$ is functional for $f$, then 
over $U$ the configuration space is a bunch of copies of the graph
 of $f$.
If the configuration space is just one copy of the graph
 of $f$ it is strongly functional.

\commentout{
A linkage $\ul$ is \fdf{functional} for a map $f\colon {\Bbb C}^n\to {\Bbb C}^m$ 
 with \rd\ $U$ if there are 
 vertices $w_1,\ldots ,w_n$ and $v_1,\ldots ,v_m$ of $\ul$ so that
 \item{1)} $\ul$ is \qf\ for $f$ and
 \item{2)} The restriction $q|\colon q^{-1}(U)\to U$ is an analytically trivial
 covering map and $U\subset q(\conf(\ul))$.
 We call $w_1,\ldots ,w_n$ the input vertices and call
 $v_1,\ldots ,v_m$ the output vertices.
 We require all input vertices to be distinct,
 and all output vertices to be distinct.
 We will, however, allow an input vertex to also be an output vertex.
 A linkage $\ul$ is \fdf{strongly functional} for a map $f\colon {\Bbb C}^n\to {\Bbb C}^m$ 
 with \rd\ $U$ if there are 
 vertices $w_1,\ldots ,w_n$ and $v_1,\ldots ,v_m$ of $\ul$ so that
 \item{1)}  $\ul$ is \qf\ for $f$ and
 \item{2)} The restriction $q|\colon q^{-1}(U)\to U$ is an analytic isomorphism.
}

It is very hard for a classical linkage to be strongly functional,
because their configuration spaces are real algebraic sets and
polynomial maps on real algebraic sets have a mod 2 degree theory
\cite{AK}.
If $\ul$ is a classical linkage and $\conf(\ul)$ is compact,
the mod 2 degree of the map $q$ must be even (since for any point $z$
outside the image of $q$, the fiber $q^{-1}(z)$ has an even number of points).
Thus if $U$ is open, the degree of the covering $q|\colon q^{-1}(U)\to U$
must be even.
So for a classical linkage to be strongly functional
with open \rd, we must have
$\conf(\ul)$ noncompact, i.e., some component of $L$ has no fixed
vertices.
A moment's thought will then confirm that only 
coordinate projections,
constant maps,
or products of the two could have strongly functional classical linkages
with open \rd.

 But \rl s are a different matter.
 Since their configuration spaces are only \qas s,
 there is no degree theory to get in the way of having strongly functional
 linkages.
 Indeed we will show the following:
 
 \proclaim{Theorem 4}. 
 (\KM)
 For any real polynomial map $g(w,\overline w)\colon {\Bbb C}^n\to {\Bbb C}^m$
 and any compact $U\subset {\Bbb C}^n$ there is a functional classical linkage
 $\ul$ for $g$ with \rd\ $U$.
 
\proclaim{Theorem 5}. 
 (King)
 For any real polynomial map $g(w,\overline w)\colon {\Bbb C}^n\to {\Bbb C}^m$
 and any compact $U\subset {\Bbb C}^n$ there is a strongly functional \rl\ 
 $\ul$ for $g$ with \rd\ $U$.

\sec{Proofs of Theorems 1, 2, and 3}
 
Given Theorems 4 and 5, we may easily prove Theorems 1, 2 and 3.
To prove Theorem 1, take any compact real algebraic set $X\subset {\Bbb C}^n$,
find a polynomial $g$ so $X=g^{-1}(0)$, then use Theorem 4 to get a
classical linkage $\ul'$ which is functional for $g$ with \rd\ $X$.
Now construct the linkage $\ul$ by taking $\ul'$ and fixing all
the output vertices to 0.
Then by Lemma 3.4, 
$$\eqalign{
\conf(\ul)&=\{ \varphi \in \conf(\ul') \mid \varphi (v_i)=0, i=1,\ldots ,m \}\cr
&=p^{-1}(0)=(f\compose q)^{-1}(0)=q^{-1}(X)\cr}$$
But $q^{-1}(X)$ is an analytically trivial covering of $X$.

%To prove Theorem 2, we will need a linkage $\ul'$ with some vertex
%$w$ which traces a straight line,
%that is, $\rh{'}w(\conf(\ul'))=[a,b]$.
%Later we will construct such a linkage explicitly
% (Watt tried in vain to find such a linkage to use for his steam engine).
%However, we may derive such a linkage from Theorem 4 as follows.
%Let $\ul''$ be a classical linkage which is functional for
%$z\to (a+b)/2 +(a-b)(z+\overline z)/4$ with \rd\ the disc $\{ |z|\le 2 \}$
%and let $w$ be its output vertex.
%Let $\ul'$ be obtained from this linkage by adding an edge
%of length 1, attaching one of its vertices to the input vertex of 
%$\ul''$ and fixing the other of its vertices to 0.
%Then the output vertex $w$ will be in the image of the unit circle under
%the above map.
%This image is the interval $[a,b]$.
%
%Now given any polynomial map $\alpha \colon [a,b]\to {\Bbb C}$ we may find a
%functional linkage $\ul'''$ for $\alpha $ with \rd\ $[a,b]$.
%Form a linkage $\ul$ from $\ul'$ and $\ul'''$
%by attaching the output vertex $w$ of $\ul'$
%to the input vertex of $\ul'''$.
%Let $v$ be the output vertex of $\ul'''$.
%Then $\rh{}v(\conf(\ul))=\alpha ([a,b])$.
%
%Of course we could have also done this more simply by 
%taking a linkage which is functional for $z\mapsto \alpha (z+\overline z)$
%and attaching to its input an edge of length $(b-a)/4$
%with a vertex fixed at $(b+a)/4$.

To prove Theorem 2,
take a linkage $\ul'$ which is functional for $z\mapsto \alpha (z+\overline z)$
and with \rd\ the circle $C$ with radius $(b-a)/4$ and center $(b+a)/4$.
Form $\ul$ from $\ul'$ by attaching to its input an edge of length $(b-a)/4$
with a vertex fixed at $(b+a)/4$.
The end of this edge will trace out the circle $C$,
but the function $z\mapsto z+\overline z$ maps $C$ onto the interval
$[a,b]$, hence the output vertex of $\ul$ will trace out
$\alpha ([a,b])$.
Alternatively, we show below that there is a linkage $\ul''$ with a vertex
$w$ which traces any interval,
that is, $\rh{''}w(\conf(\ul''))=[a,b]$.
We could hook this linkage up to a functional linkage for $\alpha $.

Theorem 3 is proven similarly to Theorem 1, using 
\rl s and Theorem 5 instead of Theorem 4.
In particular, take any compact \qas\ $X\subset {\Bbb C}^n$,
find a polynomial $g(z,\overline z)\colon {\Bbb C}^n\to {\Bbb C}^m$ 
so the image of $g$ is contained in ${\Bbb R}^m$ and
$$X=\{ x\in {\Bbb C}^n \mid g_i(x)=0, i\le k \mand g_i(x)\ge 0, i>k \}$$
Use Theorem 5 to get a \rl\
 $\ul'$ which is functional for $g$ with \rd\ $X$.
 By compactness of $X$ we may choose a real number $b$ so that
 $g_i(x)\le 2b$ for all $i>k$ and $x\in X$.

Now construct the linkage $\ul$ by taking $\ul'$, fixing all
the output vertices $v_i$ for $i\le k$ to 0,
and adding a new vertex $u_i$ 
and a new flexible edge from $v_i$ to $u_i$ for each $k<i\le m$.
We fix the new vertices $u_i$ to the point $b$ and make the
length of each new flexible edge $b$.
(We call this operation tethering $v_i$ to $u_i$.)
Since each $g_i$ is real valued, the effect of each new flexible edge is
to restrict $0\le g_i(x)\le 2b$.
Then by Lemma 3.4 and Lemma 3.5,
$$\conf(\ul)=\{ \varphi \in \conf(\ul') \mid \varphi (v_i)=0, i\le k \mand 0\le \varphi (v_i)\le 2b, k<i\le m \}
=q^{-1}(X)$$
But $q|:q^{-1}(X)\to X$ is an analytic isomorphism,
so $\conf(\ul)$ is analytically isomorphic to $X$.

To see Corollary 3.1, note that if $X$ is a compact \qas,
Theorem 3 gives a \rl\ $\ul'$ so that $\conf(\ul')$ 
is isomorphic to $X$.
Now let $\ul$ be the union of $\ul'$ with $k$
disjoint vertices.  Then $\conf(\ul)=\conf(\ul')\times {\Bbb R}^{2k}$.
On the other hand, if $\ul$ is any \rl, let
$\ul'$ be a linkage obtained by fixing one vertex
in every component of $\ul$ which does not already have
a fixed vertex.  Let $k$ be the number of such vertices
we fixed.
Then $\conf(\ul)=\conf(\ul')\times {\Bbb R}^{2k}$, and $\conf(\ul')$
is a compact \qas.

\sec{Constructing Polynomial Functional Linkages}

So it remains to prove Theorems 4 and 5.
Notice first that any polynomial map
$f(z,\overline z)\colon {\Bbb C}^n\to {\Bbb C}^m$ can be written as a
composition of the following two types of polynomial maps:
\item{$\bullet$} $g\colon {\Bbb C}^k\to {\Bbb C}^k\times {\Bbb C}$ given by
$g(z)=(z,h(z))$ where $h(z)$ is $z_i+z_j$, $z_iz_j$, $\overline z_i$,
or a constant.
\item{$\bullet$} a projection ${\Bbb C}^k\to {\Bbb C}^m $ onto some of
the coordinates, $m<k$.

This corresponds to how we calculate polynomials, as a series of
elementary arithmetic operations on partial results, and throwing
away those results no longer needed.

Consequently, because of the following lemmas,
it suffices to find 
functional linkages for the polynomials
$z+w$, $zw$, $\overline z$, a constant, and coordinate projection,
all with arbitrarily large compact \rd.
The last two are trivial. 
For the constant, just take a fixed output vertex.
For the projection, use a linkage with $k$ vertices and no edges,
make all $k$  vertices inputs, and select $\ell $ output vertices.

\proclaim{Lemma 6.1}. 
Let $\ul$ and $\ul'$ be (strong) functional linkages
for functions $f\colon {\Bbb C}^n\to {\Bbb C}^m$ and 
$g\colon {\Bbb C}^m\to {\Bbb C}^k$
with \rd s $U$ and $U'$.
Suppose that $U\cap f^{-1}(U')$ is nonempty.
Form a linkage $\ul''$ by taking the disjoint union
of $\ul$ and $\ul'$ and then identifying each
output vertex of $\ul$ with the corresponding
input vertex of $\ul'$.
Then $\ul''$ is a (strong) functional linkage
for $g\compose f$ with \rd\ $U\cap f^{-1}(U')$.

\proof
\ifdupv{
First assume that there are no duplications of the 
input vertices of $\ul'$ and also no duplications
of the output vertices of $\ul$.
}
Note that $\ul''$ is the union of 
$\ul$ and $\ul'$, and their intersection is 
a linkage with $m$ vertices and no edges.
Let $\rho _1$ and $\rho _2$ be the input and output
maps of $\ul$ and let $\rho _3$ and $\rho _4$ be the
input and output maps of $\ul'$.
By Lemma 3.2, we know that $\conf(\ul'')$ is the fiber product of
$\rho _2$ and $\rho _3$,
$$\conf(\ul'')=\{(\varphi ,\varphi ')\in \conf(\ul)\times \conf(\ul') \mid \rho _2(\varphi )=\rho _3(\varphi ')\}$$
So that $\rh{''}{\ul}$ and $\rh{''}{\ul'}$ are induced by projection.
Note that 
$$g\compose f\compose \rho _1\compose \rh{''}\ul
=g\compose \rho _2\compose \rh{''}\ul
=g\compose \rho _3\compose \rh{''}{\ul'}=\rho _4\compose \rh{''}{\ul'}$$
so $\ul''$ is \qf\ for $g\compose f$.

Now let us see that we can take the \rd\ to be $U\cap f^{-1}(U')$.
The restriction of $\rho _1$ to $U\cap f^{-1}(U')$ is an analytically
trivial cover since the restriction to $U$ is, so we only need
show that $\rh{''}\ul$ restricts to an analytically trivial
cover of $\rho _1^{-1}(U\cap f^{-1}(U'))=\rho _1^{-1}(U)\cap \rho _2^{-1}(U')$.
We know that there is a finite set $F$ and an analytic isomorphism
$\alpha \colon U'\times F\to \rho _3^{-1}(U')$ 
so that $\rho _3\compose \alpha $ 
is projection to $U'$.
Now 
$$\eqalign{
\rh{''}\ul^{-1}(\rho _2^{-1}(U'))&=\{(\varphi ,\varphi ')\mid \rho _2(\varphi )
=\rho _3(\varphi ')\in U' \}\cr
&=\{(\varphi ,\alpha (\rho _2(\varphi ),f))\mid \rho _2(\varphi )\in U' \mand f\in F\}\cr}$$
So we have an analytic trivialization 
$\alpha '\colon \rho _2^{-1}(U')\times F\to \rh{''}\ul^{-1}(\rho _2^{-1}(U'))$
given by $\alpha '(\varphi ,f)=(\varphi ,\alpha (\rho _2\varphi ,f))$.

In the strong case, note that all the covers are one-fold and hence
are analytic isomorphisms.
\ifnotdupv{\qed}

\ifdupv{
If there are duplications in the input and output vertices things get
more complicated.
For the purposes of this paper we never need to use duplicated vertices,
but for the sake of generality we provide the proof.
Let $\Delta _{ij}=\{ (z_1,\ldots ,z_m)\in {\Bbb C}^m \mid z_i=z_j\}$.
Let $v_1,\ldots ,v_m$ be the input vertices of $\ul'$ and let
$w_1,\ldots ,w_m$ be the output vertices of $\ul$.
}

\ifdupv{
Suppose $v_i=v_j$.  Then we must have $U'\subset \Delta _{ij}$.
Also, in $\ul''$ we end up identifying $w_i$ with $w_j$.
Let us first see whether we can do so.
Suppose $w$ is another vertex so that $ww_i$ and $ww_j$
are both edges of $\ul$.  
Since $U\cap f^{-1}(U')$ is nonempty,
there is a $\varphi \in \conf(\ul)$ so that $\rho _1(\varphi )\in U\cap f^{-1}(U')$.
Hence $\rho _2(\varphi )=f(\rho _1(\varphi ))\in U'\subset \Delta _{ij}$, and so $\varphi (w_i)=\varphi (w_j)$.
So $\ell (ww_i)=|\varphi (w)-\varphi (w_i)|=|\varphi (w)-\varphi (w_j)|=\ell (ww_j)$.
So in $\ul''$ we may identify the edges $ww_i$ and $ww_j$
since they have the same length.
There could not be an edge $w_iw_j$ since 
$0\neq \ell (w_iw_j)=|\varphi (w_i)-\varphi (w_j)|=0$.
So by Lemma 3.7 we are allowed to take the quotient linkage
$\ul_1$ of $\ul$, identifying $w_i$ and $w_j$.
By Lemma 3.7 we also see that $\ul_1$ is still functional for 
$f$ but the domain has been restricted to $U\cap f^{-1}(\Delta _{ij})$.
Do this identification for each pair $i,j$ with $v_i=v_j$ and we
eventually get a functional linkage $\ul_2$ for $f$ with domain
$U_2=U\cap f^{-1}(\Delta )$ for some $\Delta \supset U'$.
}

\ifdupv{
Now suppose $w_i=w_j$.
Then we must have $f(U)\subset \Delta _{ij}$.
Also, in $\ul''$ we end up identifying $v_i$ with $v_j$.
Let us see whether we can do so.
Suppose $v$ is another vertex so that $vv_i$ and $vv_j$
are both edges of $\ul'$.  
Since $U\cap f^{-1}(U')$ is nonempty, we know $\Delta _{ij}\cap U'$
is nonempty, so
there is a $\varphi \in \conf(\ul')$ so that $\rho _3(\varphi )\in \Delta _{ij}\cap U'$,
and hence $\varphi (v_i)=\varphi (v_j)$.
So as above, Lemma 3.7 will allow us to take the quotient linkage
$\ul'_1$ identifying $v_i$ and $v_j$.
By Lemma 3.7 we also see that $\ul'_1$ is still functional for 
$g$ but the domain has been restricted to $U'\cap \Delta _{ij}$.
Do this identification for each pair $i,j$ with $w_i=w_j$ and we
eventually get a functional linkage $\ul'_2$ for $g$ with domain
$U'_2=U'\cap \Delta '$ for some $\Delta '\supset f(U)$.
}

\ifdupv{
After doing all this, we have $\ul''$ is the union of $\ul_2$ and $\ul'_2$,
and we may finish the proof as above.
The only thing to check is that $U_2\cap f^{-1}(U'_2)=U\cap f^{-1}(U')$.
But $U_2\cap f^{-1}(U'_2)=U\cap f^{-1}(\Delta )\cap f^{-1}(U')\cap f^{-1}(\Delta ')=U\cap f^{-1}(U')$
since $U'\subset \Delta $ and $U\subset f^{-1}(\Delta ')$.
\qed}

\proclaim{Lemma 6.2}. 
Let $\ul$ be a (strong) functional linkage
for a function $f\colon {\Bbb C}^n\to {\Bbb C}^m$
with \rd\ $U$.
\ifnotdupv{
Suppose that no vertex of $\ul$ is both an input and an
output vertex.}
Form a functional linkage $\ul'$ by taking  $\ul$ 
 but taking the output vertices of $\ul'$
to be the concatination of the input and output vertices of $\ul$.
Then $\ul'$ is a (strong) functional linkage
for $(f,g)\colon {\Bbb C}^n\to {\Bbb C}^m\times {\Bbb C}^n$ with \rd\ $U$,
where $g$ is the identity.

This lemma is trivial.
If we wished we could generalize it to a way of combining
functional linkages for any $f$ and $g$ to a linkage
for $(f,g)$, just identify the inputs.
But we have no need to do so.

\sec{Elementary Polynomial Functional Linkages}

So it suffices to find 
functional linkages for the polynomials
$z+w$, $zw$, and $\overline z$,
 with arbitrarily large compact \rd.
In fact we will not directly construct functional linkages 
for $z+w$, $zw$, and $\overline z$, but will do so for
other functions which can be composed to give them.
In particular we will construct functional linkages for:
\medskip
\item{1)} Translation:  $z\mapsto z+z_0$, with \rd\ any compactum in ${\Bbb C}$.

\item{2)} Real scalar multiplication: $z\mapsto \lambda z$, with \rd\ a disc
  $\{ |z-z_0|\le r \}$ for some $z_0$ and $r$ as large as we wish.
  
 \item{3)} Average:  $(z,w)\mapsto (z+w)/2$, with \rd\
 $\{(z,w)\mid |z-z_0|\le r, |w+z_0|\le r\}$  for some $z_0$ and $r$ as large as we wish.
 
 \item{4)} Inversion through a circle:  $z\mapsto t^2/\overline z$, with \rd\
 $\{ |z-z_0|\le r \}$ for any specified  $z_0$ with $|z_0|=t$ and any $r\le t/2$.
 
 \item{5)} Complex conjugation: $z\mapsto \overline z$, with \rd\
$\{ |z-z_0|\le r \}$ for some $z_0$ and $r$ as large as we wish.

\medbreak
From these functional linkages, we may compose to get functional linkages for
$z+w$, $zw$, and $\overline z$ with any compact \rd\ $K$.
Note first that we may always restrict the domain of a
functional linkage further, so it suffices to find functional linkages with
arbitrarily large compact \rd s, for example
(products of) balls of radius $r$.
\medskip
\item{2.1)} To get $z\mapsto \lambda z$ with \rd\ $|z|\le r$ for $\lambda $ real, first
use 2) above to get a functional linkage for $z\mapsto \lambda z$, with \rd\
  $\{ |z-z_0|\le r \}$.  Then use 1) to get a functional linkage for
  $z\mapsto z+z_0$ with \rd\ $|z|\le r$.
  Using Lemma 6.1, compose these two to get a functional
  linkage for $z\mapsto \lambda z+\lambda z_0$ with \rd\ $|z|\le r$.
  Now using 1) and Lemma 6.1, compose with a translation by $-\lambda z_0$
  to get our desired linkage for $z\mapsto \lambda z$ with \rd\ $|z|\le r$.

\item{3.1)} To get $(z,w)\mapsto z+w$ with \rd\ 
$|z|\le r$, $|w|\le r$, find a functional
linkage for the average 3) above.  Then using 1), find functional linkages
for $z\mapsto z+z_0$ and $z\mapsto z-z_0$, both with \rd\ $|z|\le r$.
Their disjoint union is functional for
$(z,w)\mapsto (z+z_0,w-z_0)$ with \rd\ $|z|\le r$, $|w|\le r$.
Using Lemma 6.1 and composing with the first linkage, we get
a functional linkage for $(z,w)\mapsto (z+w)/2$ with \rd\ $|z|\le r$, $|w|\le r$.
Now using 2.1 with $\lambda =2$ and Lemma 6.1, we get a functional
linkage for $(z,w)\mapsto z+w$ with \rd\ $|z|\le r$, $|w|\le r$.

\item{5.1)} To get $z\mapsto \overline z$, with \rd\ $|z|\le r$,
first use  5) to get a functional linkage for $z\mapsto \overline z$
with \rd\ $|z-z_0|\le r$.
Then compose with translations by $z_0$ and $-\overline z_0$.

\item{6)}  We will do some algebraic manipulation to get multiplication
$(z,w)\mapsto zw$.
First, note that $zw=((z+w)^2-(z-w)^2)/2$.
So it suffices to find a functional linkage for $z\mapsto z^2$ with
\rd\ $|z|\le r$.
Next, note that if $h(z)=t^2/\overline z$ then
$$t^2-th((h(t+z)+h(t-z))/2)=z^2$$
Consequently, we may get a functional linkage for $z^2$ by 
using Lemma 6.1 and functional linkages we have found
above.
In particular, the reader who wishes to work out the details
will find it is good to use (3 copies of) a functional linkage for
$h$ by choosing some $t>2r$, and specifying that its \rd\ be
$\{|z-t|\le r\}$.

\medbreak
So we have now reduced to finding functional linkages 1)-5) above.
In doing so, the following Lemma will be useful.
Its proof may be safely left to the reader.
It is, for example, a special case of the theorem that
a proper submersion is a locally trivial fibration.

\proclaim{Lemma 6.3}. 
Let $f\colon M\to {\Bbb R}^m$ be a smooth map from a compact
$m$ dimensional manifold with boundary.
Let $S\subset M$ be the set of critical points of $f$,
the points where $df$ has rank $<m$.
Let $U$ be any connected component of ${\Bbb R}^m-f(S\cup \partial M)$.
Then $f|\colon f^{-1}(U)\to U$ is a covering projection.

In our usage, $m=2$, $f$ is analytic, and $U$
is often contractible, so $f$ restricts to an analytically trivial covering
of $U$, thus $f^{-1}(U)$ is analytically isomorphic to $U\times $ a finite set.
As another application, we will use the consequence that
$f(M)$ is the union of $f(S\cup \partial M)$ and some connected components
of ${\Bbb R}^m-f(S\cup \partial M)$.

\sec{A simple Linkage, a key to understanding more complicated Linkages}

It will be useful to look first at a very simple linkage $\ul$,
as shown in the left half of Figure 6.
The vertex $A$ is fixed at some point $z_1$,
but $B$ and $C$ are not fixed.
The two edges $AB$ and $BC$ are rigid.
We assume that $b\le a$.

\midinsert\vskip 1.5in
\includegraphics{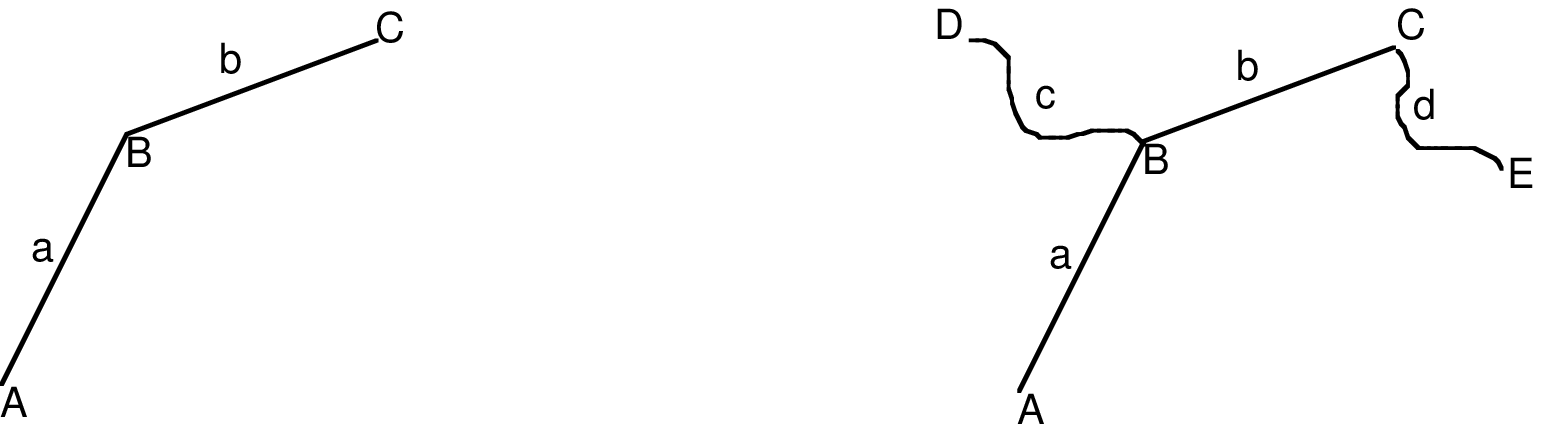}
\capshun 6 { The Tensor Lamp Linkage}
\endinsert

Clearly $\conf(\ul)$ is a torus $S^1\times S^1$,
we may identify $(u,v)\in S^1\times S^1$ with
 $\varphi _{uv}$ where $\varphi _{uv}(A)=z_1$, $\varphi _{uv}(B)=z_1+au$,
 and $\varphi _{uv}(C)=z_1+au+bv$.
Note that $\rh{}C\colon \conf(\ul)\to {\Bbb C}$ is then the map
$\rh{}C(u,v)=z_1+au+bv$ which has critical set
where $u=\pm v$.
The image of the critical set is the two circles
$\{ |z-z_1|=a-b \}\cup \{ |z-z_1|=a+b \}$,
so by Lemma 6.3 we see that the image of
$\rh{}C$ is the annulus  $\{ a-b\le |z-z_1|\le a+b\}$
and moreover $\rh{}C$ restricts to a double cover
of the open annulus $\{ a-b<|z-z_1|<a+b \}$.
In fact this double cover is analytically trivial.
(The masochistic reader may verify that the inverse of
$\rh{}C$ takes $B$ to the points
$z_1+(z-z_1)(\alpha +a^2-b^2\pm \beta )/(2\alpha )$ where
$\alpha =|z-z_1|^2$ and $\beta =\sqrt{(\alpha +a^2-b^2)^2-4a^2\alpha }$.)
In applications below, we will usually only focus on some disc
inside the annulus $\{|z-z_0|\le r\}$ where, say $z_0=z_1+a$
and $0<r<b$.
Then $\rh{}C$ restricts to an analytically trivial
double cover of this disc.

When working with \rl s, it will be convenient to modify
this linkage so that $\rh{}C$ is an analytic isomorphism
to some disc $\{|z-z_0|\le d\}$.
We do this by tethering the vertices $B$ and $C$
 to fixed vertices $D$ and $E$ so that their movement
is restricted, see the \rl\ on the right half of Figure 6.
Consider first the sublinkage $\ul'$ formed by $A$, $B$, $C$,
and $D$, with rigid edges $AB$ and $BC$, and a flexible
edge $BD$ of length $c$, where $D$ is fixed at some point $z_2$
(and $A$ is fixed at $z_1$ as before).
We have 
$$\conf(\ul')=\rh{}B^{-1}(\{|z-z_2|\le c\})$$
$$=\{(u,v)\in S^1\times S^1 \mid |z_1+au-z_2|\le c \}=T\times S^1$$
for some arc $T$ of $S^1$ as long as we choose $z_2$
and $c$ so that $c-a<|z_1-z_2|<c+a$.
For convenience, we choose $c=\sqrt 2 a$ and $z_2=z_1+aw_0$
for some $w_0\in S^1$.
Then $T$ will be the semicircle between $\pm \sqrt{-1}w_0$
which contains $w_0$.
By Lemma 6.3,
we know that $\rh{}C$ restricts to an analytically trivial
covering of $\{|z-z_1-\sqrt{-1}aw_0|<b\}$.
But by checking the inverse image of a point, for example
$z_1+\sqrt{-1}aw_0$, we see that it is a one-fold cover,
hence an analytic isomorphism.
So now in $\ul$, if we fix $E$ at $z_1+\sqrt{-1}aw_0$
and pick $d<b$, we see that
$\rh{}C\colon \conf(\ul)\to {\Bbb C}$ is an analytic isomorphism to its
image $|z-z_1-\sqrt{-1}aw_0|\le d$.

\sec{A Functional Linkage for Translation }

Now let us find a functional linkage for translation 1) above
with \rd\ $|z|\le r$.
Consider the linkages $\ul$ in Figure 7, which we will show to be functional
for $z\mapsto z+z_0$ with \rd\ $|z|\le r$.
The right hand \rl\ will be strongly functional.

\midinsert\vskip 1.6in
\includegraphics{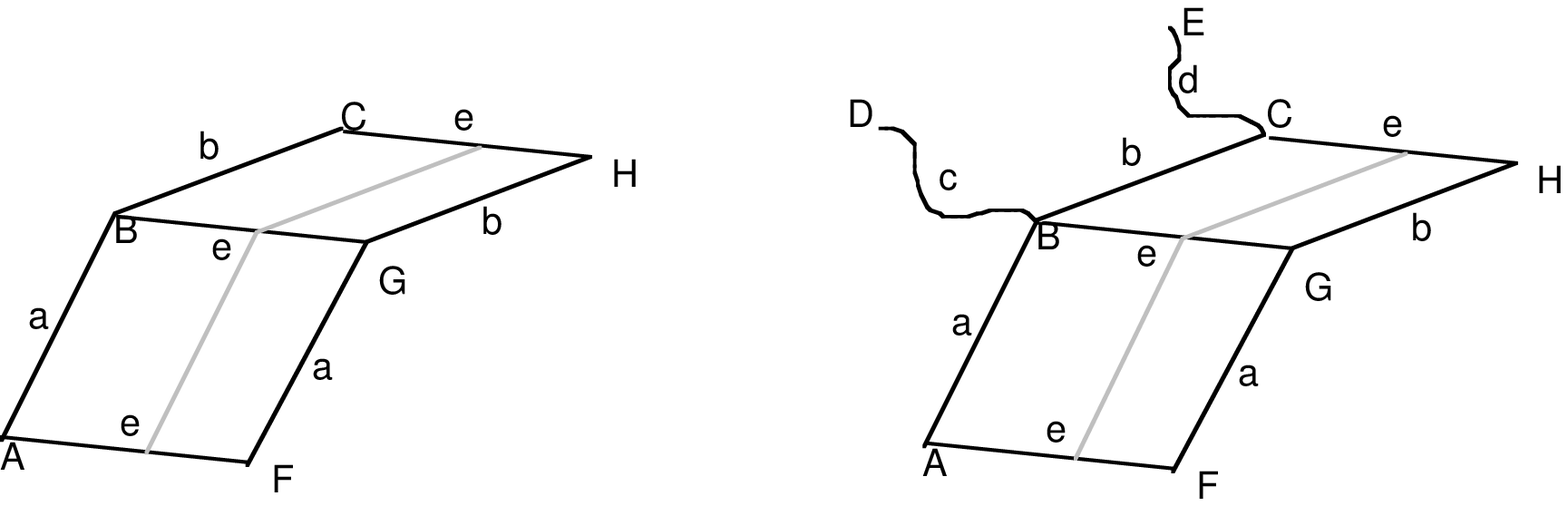}
\capshun 7 {Translation Linkage }
\endinsert

Choose $a>2r$, and let $e=|z_0|$.
The vertices $A$ and $F$ are fixed at $z_1$ and $z_1+z_0$
respectively, where $z_1$ will be determined later.
We let $C$ be the input vertex and $H$ be the output vertex.
The parallelograms $ABGF$ and $BCHG$ are rigidified.
Thus $AF$, $BG$, and $CH$ are parallel, and so
for any $\varphi \in \conf(\ul)$ we must have $\varphi (H)=\varphi (C)+z_0$.
So $\ul$ is \qf\ and we must only check that the \rd\ can be $|z|\le r$.
Note that the linkages $\ul'$ of Figure 6 are sublinkages of $\ul$.
We claim by Lemma 3.6 
that $\rh{}{\ul'}\colon \conf(\ul)\to \conf(\ul')$ is an isomorphism.
This is because the positions of $F$, $G$, $H$, and the other three
unnamed vertices used to rigidify the quadrilaterals are all
polynomial functions of the positions of $A$, $B$, and $C$.
Now the fact that $\rh{}C$ doubly covers $|z|\le r$ (for the left hand
classical linkage) or singly covers $|z|\le r$ (for the right hand
\rl) follows from the above discussion of $\conf(\ul')$
in Figure 6,
as long as we make appropriate choices of $z_1$, $a$, $b$,
and $w_0$.
For example, we may choose $b$ so $r<b<a$,
choose $w_0=1$ and $z_1=-a\sqrt{-1}$.

\sec{A Functional Linkage for real scalar Multiplication}

Now let us find a functional linkage $\ul$ for real scalar multiplication
2) above.
We have already seen this linkage, the pantograph in Figure 5.
We consider 3 cases, $\lambda >1$, $0<\lambda <1$, and $\lambda <0$.
The remaining cases $\lambda =0$ or $\lambda =1$ are trivial functions which 
have trivial functional linkages.

If $\lambda >1$ we take $c=\lambda -1$, let $B$ be the input vertex and let
$C$ be the output vertex, and fix $A$ at 0.
Note that for any $\varphi \in \conf(\ul)$ we have $\varphi (C)=\lambda \varphi (B)$.
So $\ul$ is \qf.
Consider the sublinkage $\ul'$, with vertices $A$, $E$, and $B$
and edges $AE$ and $EB$.
Note by Lemma 3.6
that $\rh{}{\ul'}\colon \conf(\ul)\to \conf(\ul')$ is an isomorphism
since the positions of $D$, $F$, and $C$ are polynomial functions
of $A$, $E$, and $B$.
By the discussion of the linkage in Figure 6,
we know that if $a$ and $b$ are chosen appropriately,
then $\rh{}B$ double covers some disc $|z-z_0|\le r$.
To get a strong functional linkage, we add two vertices and
tether $E$ and $B$ to them with appropriate length cables.
By the discussion of the right hand linkage of Figure 6,
we know that $\rh{}B$ singly covers some disc $|z-z_0|\le r$.

If $0<\lambda <1$, we take $c=1/\lambda -1$, let $C$ be the input vertex,
let $B$ be the output vertex, and fix $A$ at 0.
By considering the sublinkage with vertices $A$, $D$, and $F$,
we see as above that with appropriate choices of $a$ and $b$,
the linkage will be functional for $z\mapsto \lambda z$ with \rd\
$|z-z_0|\le r$.
To get a strongly functional linkage, tether $D$ and $C$
appropriately.

If $\lambda <0$, we take $c=-\lambda $, let $A$ be the input vertex, 
$C$ be the output vertex, and fix $B$ at 0.
Letting $\ul'$ be the sublinkage with vertices $A$, $E$, and $B$,
we see as above that $\ul$ is functional for $z\mapsto \lambda z$ with \rd\
$|z-z_0|\le r$.
To get a strongly functional linkage, tether $E$ and $A$
appropriately.

\sec{A Functional Linkage for the Average}

Now let us find a functional linkage $\ul$ for the average 3).
Again $\ul$ will be the pantograph of Figure 5.
The input vertices will be $A$ and $C$.
The output vertex will be $B$.  
We let $c=1$.
Note $\ul$ is \qf\ for $(z,w)\mapsto (z+w)/2$.

Let $\ul'$ be the linkage obtained from $\ul$ by
fixing $A$ to the point $0$.
We know by Lemma 3.3 that there is an isomorphism
$\alpha \colon \conf(\ul')\times {\Bbb C}\to \conf(\ul)$
where $\alpha (\varphi ,w)(v)=\varphi (v)+w$ for any vertex $v$ of $\ul$.
Thus $\rh{}C\compose \alpha (\varphi ,w)=w+\rh{'}C(\varphi )$ and 
$\rh{}A\compose \alpha (\varphi ,w)=w$.
Note that $\ul'$ is the functional linkage for $z\mapsto z/2$
which we have already studied.
So we know for appropriate choices of $a$ and $b$, that
$\rh{'}C\colon \conf(\ul')\to {\Bbb C}$ double covers
 some disc $|z-2z_0|\le 2r$.
Consequently, $(\rh{}C,\rh{}A)\colon \conf(\ul)\to {\Bbb C}^2$ double covers
 the set $|z-2z_0-w|\le 2r$,
since $(\rh{}C,\rh{}A)\compose \alpha $ does.
In particular, $(\rh{}C,\rh{}A)$ is an analytically trivial
double cover of the subset $|z-z_0|\le r$, $|w+z_0|\le r$.

Now let us construct a strong functional \rl\ for the average.
As usual, we will do so by tethering some vertices of $\ul$,
in particular, $A$, $C$, and $D$.
In fact $|z_0|>r$ already, but in any case
we may shrink $r$ so this is true.
We will first tether $A$ to a fixed vertex at $-z_0$ with a
cable of length $r$.
Then tether $C$ to a fixed vertex at $z_0$ with a
cable of length $r$.
From the above analysis, the resulting linkage $\ul'$
will have $\conf(\ul')$ be an analytically trivial double
cover of the product of discs $K=\{(z,w)\mid |z-z_0|\le r, |w+z_0|\le r\}$.
We can compute the inverse of this double cover explicitly,
for convenience we choose our linkage so that $a=b$,
and of course $c=1$.
If $(\rh{}C,\rh{}A)(\varphi )=(z,w)$ for $(z,w)\in K$,
then 
$$\varphi (D)=\gamma (z,w,f)=(z+w)/2+f\sqrt{-1}\beta (z,w)(w-z)/|w-z|$$
where $f=\pm 1$ and 
where $\beta \colon K\to (0,\infty )$ is 
$$\beta (z,w)=\sqrt{4a^2-|w-z|^2/4}$$
Note $\gamma (z_0,-z_0,1)\neq \gamma (z_0,-z_0,-1)$.
So by shrinking $r$ if necessary we may find a $d$ so that
for all $(z,w)\in K$ we have $|\gamma (z,w,1)-\gamma (z_0,-z_0,1)|<d$
and $|\gamma (z,w,-1)-\gamma (z_0,-z_0,1)|\ge d$.
So if we tether $D$ to a fixed vertex at $\gamma (z_0,-z_0,1)$
with a cable of length $d$,
the resulting \rl\ will be strongly functional for $(z,w)\mapsto (z+w)/2$
with \rd\ $K$.

The only problem is that since we shrunk $r$, the \rd\ is no longer
as large as we want.
So we do a final step, where we take our \rl\
and rescale it by some large factor $N$.
In other words, all edge lengths are multiplied by $N$,
and if a vertex $v$ is fixed at some $z$,
we instead fix it at $Nz$.
The resulting linkage is still strongly functional for
$z\mapsto (z+w)/2$, since $N(z/N+w/N)/2=(z+w)/2$
but its \rd\ is $\{(z,w)\mid |z-Nz_0|\le Nr, |w+Nz_0|\le Nr \}$
and hence as big as we wish.

\sec{A Functional Linkage for Inversion through a Circle}

We must now find a functional linkage $\ul$ for inversion 
through the circle 4).
For this we will use the \PC\ inversor of Figure 8.
The linkage at the left is the full linkage,
the one at the right just has the basics.
The extra vertices and edges are only needed to
eliminate some degenerate configurations
when $B$ and $C$ coincide.
This is slightly modified from the linkage given
in \cite{KM}, to make $\conf(\ul)$ easier to compute.

% insert Figure 8
\midinsert\vskip 1.75in
\includegraphics{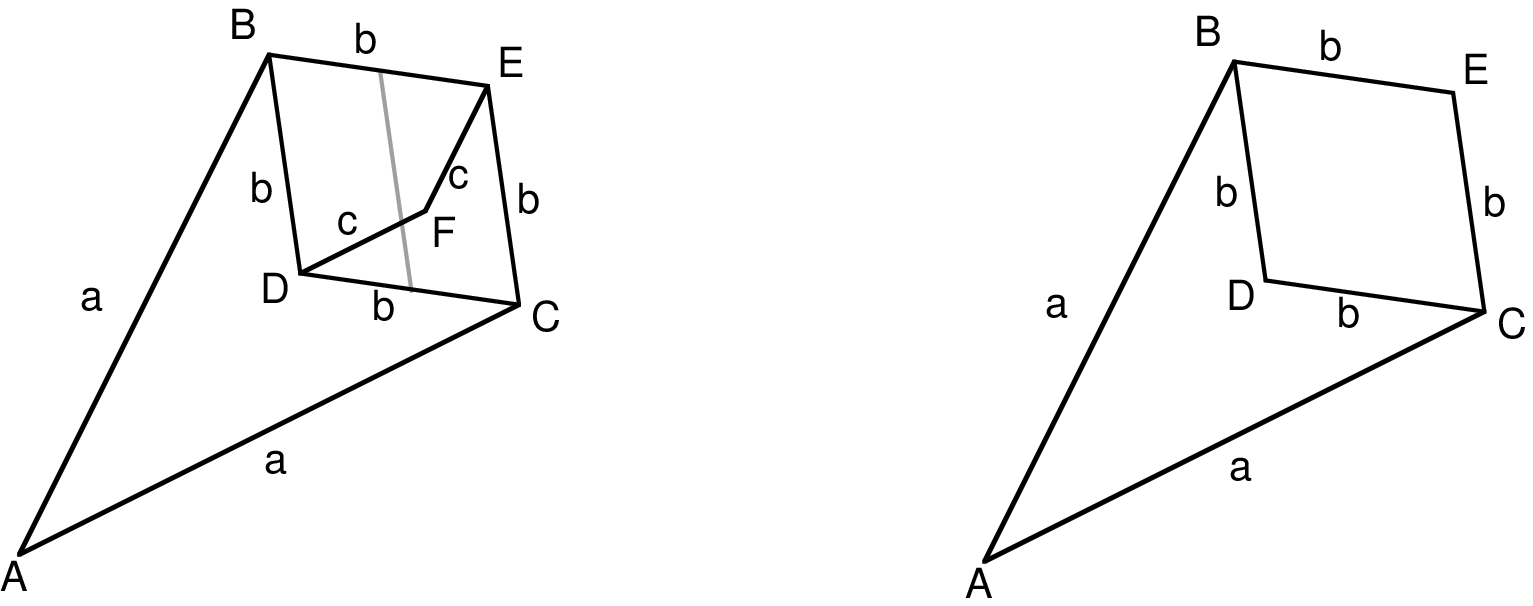}
\capshun 8 { The \PC\ Inversor}
\endinsert

We fix $A$ at 0, set $t^2=a^2-b^2$, $c<b<a$, let the input vertex be
$D$ and the output be $E$.
Let us see why $\ul$ is \qf\ for $z\mapsto t^2/\overline z=t^2z/|z|^2$.
If $\varphi \in \conf(\ul)$ note that $\varphi (D)$ is a multiple of $\varphi (E)$.
This follows from the fact that the lines from $\varphi (A)$, $\varphi (D)$, and $\varphi (E)$
to the midpoint of the line segment
$\varphi (B)\varphi (C)$ are all perpendicular to $\varphi (B)\varphi (C)$,
hence they are colinear.
Solving triangles shows that if $s=|\varphi (B)-\varphi (C)|/2$, then 
$$|\varphi (D)|=\sqrt{a^2-s^2}\pm \sqrt{b^2-s^2}$$
 $$|\varphi (E)|=\sqrt{a^2-s^2}\mp \sqrt{b^2-s^2}$$
from which we see that $|\varphi (D)|\,|\varphi (E)|=(a^2-b^2)$.
So $\ul$ is \qf\ for $z\mapsto t^2/\overline z$.

Now let us check that for any $0<r\le t/2$ then
for appropriate choices of $a$, $b$, and $c$,
we may take the \rd\ to be $|z-t|\le r$. 
By being more careful below, we would really only need $r<t$,
but it is not worth the extra bother.
As $r$ gets closer to $t$, the lengths $a$, $b$, and $c$
need to get much bigger than $t$.

Note first that for any $\varphi \in \conf(\ul)$ we have $|\varphi (F)|=\sqrt{t^2+c^2}$.
This is because, if $u$ is the distance from $(\varphi (B)+\varphi (C))/2$ to $\varphi (F)$,
then since $(\varphi (B)+\varphi (C))/2=(\varphi (D)+\varphi (E))/2$ we have
$$\eqalign{
|\varphi (F)|^2&=u^2+((|\varphi (D)|+|\varphi (E)|)/2)^2\cr
 c^2&=u^2+((|\varphi (D)|-|\varphi (E)|)/2)^2\cr}$$
Subtracting these equations gives 
$$|\varphi (F)|^2-c^2=|\varphi (D)|\,|\varphi (E)|=t^2$$

It will be convenient to consider another linkage $\ul'$ obtained
from $\ul$ by adding an edge between $A$ and $F$
with length $d=\sqrt{t^2+c^2}$.
Since $|\varphi (F)|=d$ for any $\varphi \in \conf(\ul)$,
the map $\rh{'}{\ul}\colon \conf(\ul')\to \conf(\ul)$ is an isomorphism.
Now $\ul'$ contains a sublinkage $\ul''$ with vertices
$A$, $F$, and $D$, and edges $AF$ and $FD$.
This is just our old friend from Figure 6.
However, the map $\rh{'}{\ul''}\colon \conf(\ul')\to \conf(\ul'')$ is no longer an
isomorphism as it has been in previous linkages.
We will show that it is an analytically trivial double cover.
(The nontrivial covering translation is obtained from 
reflection about the line $ADE$.)
We see this by solving for $\varphi (B)$, $\varphi (C)$, and $\varphi (E)$ in terms
of $\varphi (D)$ and $\varphi (F)$.
In particular,  if $f=\pm 1$ and  
$$\eqalign{
\beta (u,v)&=\sqrt{b^2-|u-v|^2/4}\cr
\gamma (u,v,f)&=(u+v)/2+f\sqrt{-1}\beta (u,v)u/|u|\cr}$$
then
$$\eqalign{
\varphi (E)&=t^2/\overline{\varphi (D)}\cr
\varphi (B)&=\gamma (\varphi (D),\varphi (E),f)\cr
\varphi (C)&=\varphi (D)+\varphi (E)-\varphi (B)\cr
}$$

 The function $\beta (\varphi (D),\varphi (E))$
 is analytic on $\conf(\ul'')$
since 
$$|\varphi (D)-\varphi (E)|^2/4\le (2c)^2/4<b^2$$
Hence the above maps give an analytic isomorphism
$$\alpha \colon \conf(\ul'')\times \{1,-1\}\to \conf(\ul')$$
Putting this all together with the analysis of Figure 6, we see that 
$\rh{}D$ restricts to a fourfold cover of the annulus
$d-c<|z|<d+c$.
So  pick $c=t$, then $d=\sqrt 2 t$.
If $|z_0|=t$ and $|z-z_0|\le r\le t/2$, then $d-c<t/2\le |z|\le 3t/2<c+d$.
So  $\rh{}D$ restricts to a fourfold cover of the disc $|z-z_0|\le r$,
and thus $\ul$ is functional for inversion with \rd\ $|z-z_0|\le r$.

We could get a strong functional \rl\ by taking $\ul$ above,
and tethering $F$, $D$, and $B$ to fixed vertices,
using appropriate length cables.
(Those for $F$ and $D$ could be chosen as in Figure 6,
and for example if  $b=6t$, $a=\sqrt{37}t$, and $e=5t$
we could tether $B$ to $z_0(1+b\sqrt{-1}$ with a cable of
length $e$.)

However, it is easier to modify the linkage $\ul$ by
eliminating the vertex $F$ and edges $DF$ and $EF$,
and replacing them with a flexible edge $DE$ of length
$2c$.  This has the same effect (preventing $B$ and $C$
from coinciding) but the configuration space is simpler.
So let $\ul'$ be the \rl\ obtained by the above modification,
and also by tethering $C$ and $D$ to fixed vertices.
By the analysis of Figure 6, we may do this tethering of 
$C$ and $D$ in such a way that 
 if $\ul''$ is the sublinkage formed by $AC$ and $CD$
then
  $\rh{''}D\colon \conf(\ul'')\to {\Bbb C} $
is an isomorphism to its image  $|z-z_0'|\le r'$
where $z_0'=\sqrt{-1} aw_0$ and $w_0\in S^1$ can be chosen
arbitrarily, and $r'$ is any positive number less than $b$.
Choose $w_0=-\sqrt{-1}z_0/t$, and choose for example
$a=5t/3$, $b=4t/3$, $c=t$, $r'=7t/6$.
Then if $|z-z_0|\le r$ we have 
$$|z-z_0'|\le r+|z_0-z_0'|=r+a-t<7t/6=r'$$
so $\rh{''}D$ restricts to an analytic isomorphism to 
$\{|z-z_0|\le r\}$.
In fact we may as well tether $D$ to a vertex fixed at $z_0$
by a cable of length $r$ and assume $\rh{''}D$ is an analytic isomorphism to 
$\{|z-z_0|\le r\}$.
Note that if $|z-z_0|\le r$, then 
$$|z-t^2/\overline z|=|\,|z|-t^2/|z|\,|\le t^2/(t-r)-(t-r)\le 3t/2\le 2c$$
Consequently, by Lemma 3.6
the map $\rh{'}{\ul''}$ is an analytic isomorphism,
since it is onto, and if $\varphi \in \conf(\ul')$ we know
$\varphi (E)=t^2/\varphi (D)$ and 
$\varphi (B)= \varphi (D)+\varphi (E)- \varphi (C)$.
So $\ul'$ is strongly functional for inversion with \rd\
$\{|z-z_0|\le r\}$.

\sec{How to draw a straight line}
The only remaining function, complex conjugation, 
will require that we find a linkage so that some vertex
is constrained to lie in some line segment.
Watt tried to find such a linkage for his steam engine, but ended
up getting only an approximate straight line.
But such a linkage was found in the 1860's 
by \PC\ and was actually
used for a while soon after in a ventilating scheme for the British house of
Parliament \cite{CR}.

So in this section we will construct a linkage $\ul$
with a
vertex $A$ so that, $\rh{}A(\conf(\ul))=[a,b]$, and $\rh{}A$ restricted to
$(a,b)$ is an analytically trivial cover.
In the \rl\ case we ask that $\rh{}A$ be an analytic isomorphism
from $\conf(\ul)$ to $[a,b]$.
This linkage $\ul$ is obtained by taking the input of a functional linkage for
inversion through a circle and forcing this input to lie in a circle
going through the origin.
But when we invert a circle through the origin, we get a
straight line.
Now it is just a matter of translating and rotating it and rescaling,
to make this line be any interval on the real axis.
This linkage $\ul'$ is shown in Figure 9.

% insert Fig 10
\midinsert\vskip 2in
\includegraphics{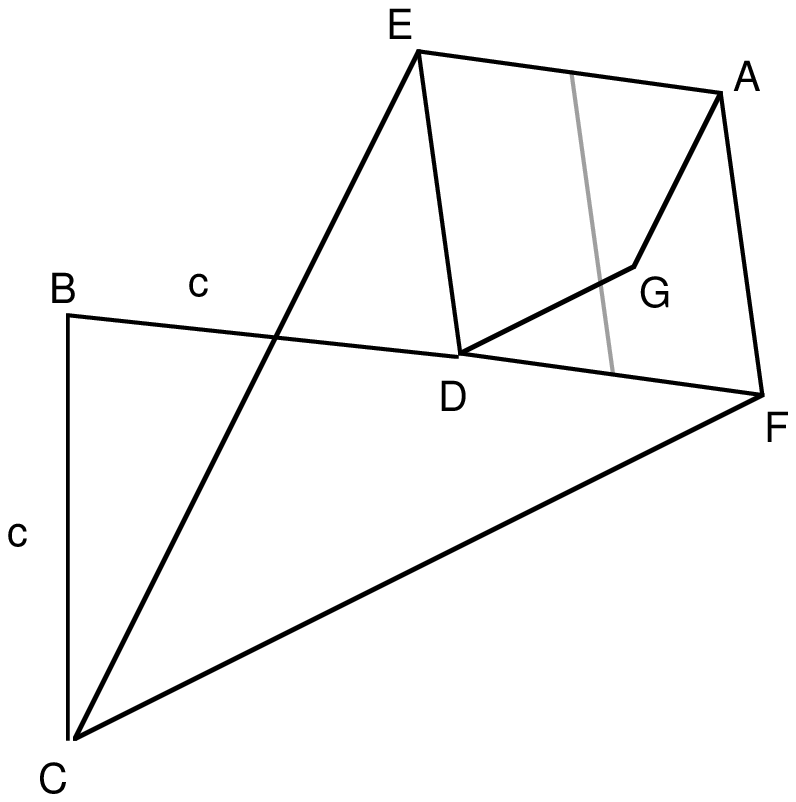}
\capshun {9} { With $B$ and $C$ fixed, $A$ will trace out a straight line}
\endinsert

This linkage is a union of two sublinkages $\ul''$ and $\ul'$.
The linkage $\ul''$ is a functional linkage for inversion through
a circle of radius $2c$ as in Figure 8,
only we translate $C$ to some point $z_0=x_0-2c\sqrt{-1}$
where the real number $x_0$ is to be determined later.
(The only effect of $x_0$ is to translate the output $A$ along the real axis).
The \rd\ of this functional linkage is $|z-x_0|\le c$.
After taking account of the translation of $C$, $\ul''$
is functional for $f(z)=z_0+4c^2/(\overline{z-z_0})$.

The linkage $\ul'$ has vertices $B$, $C$, and $D$, and
edges $BC$ and $BD$ of length $c$.
We fix $C$ at the point $z_0$ and fix $B$ at the point $z_0+c\sqrt{-1}$

Let $T$ be the circle of radius
$c$ around $z_0+c\sqrt{-1}$.
Now $\conf(\ul')=T$ and we know $\conf(\ul)$ is the fiber product of 
$\rh{'}D\colon \conf(\ul')\to {\Bbb C}$ (which is just inclusion of $T$)
 and $\rh{''}D\colon \conf(\ul'')\to {\Bbb C}$  (which is the
input map of $\ul''$).
Hence $\conf(\ul)=\rh{''}D^{-1}(T)\subset \conf(\ul'')$.

So in the classical linkage case we must do a bit more work to see what this
configuration space is.
Let $d$ be the length of $DG$ and $AG$.
Recall from the analysis of the linkage of Figure 8 that $\conf(\ul'')$
is an analytically trivial double cover of the linkage $\ul'''$
with vertices $C$, $D$, and $G$, and with edges $DG$ of length $d$
and $CG$ of length $e=\sqrt{4c^2+d^2}$.
From the analysis of Figure 6, we then see that
 the images 
$$\rh{''}D(\conf(\ul''))=\rh{'''}D(\conf(\ul'''))=\{e-d\le |z-z_0|\le e+d\}$$
Moreover, $\rh{''}D\colon \conf(\ul'')\to {\Bbb C}$ restricts to an analytically trivial
fourfold cover of the open annulus $\{e-d<|z-z_0|<e+d\}$.

As a consequence, $\rh{}A(\conf(\ul))=f(T\cap \{e-d\le |z-z_0|\le e+d\})$.
Also, since $f|\colon {\Bbb C}-z_0\to {\Bbb C}-z_0$ is a diffeomorphism and
$\rh{}A=f\compose \rh{}D$, we know that $\rh{}A$ restricts to a fourfold cover
over $f(T\cap \{e-d<|z-z_0|<e+d\})$.
But $T\cap \{e-d\le |z-z_0|\le e+d\}$ is an arc of the circle $T$,
and its image under $f$ is a line segment on the real axis.
By adjusting the constants $x_0$, $d$ and $c$ we can get any 
line segment of the real axis we desire.
(The line segment we get is the one between the
two points $x_0\pm \sqrt{2d^2+2de}$).

The analysis in the \rl\ case is a bit easier.
Choose $\ul''$ strongly functional.
By tethering $D$ to a vertex fixed at $x_0$ with a cable of length $c$,
we may then assume that $\rh{''}D\colon \conf(\ul'')\to {\Bbb C}$ 
is an analytic isomorphism to its image $\{|z-x_0|\le c\}$.
Then $\conf(\ul)=\rh{''}D^{-1}(T)$ is analytically isomorphic to
$T\cap \{|z-x_0|\le c\}$ which is again an arc of $T$,
which $f$ takes to a line segment.
(In fact to the segment $[x_0-2c/\sqrt 3,x_0+2c/\sqrt 3]$.)
So $\rh{}A$ is an analytic isomorphism from $\conf(\ul)$
to this segment.

%\vfill\eject
\sec{A Functional Linkage for Complex Conjugation}

Finally we need a functional linkage $\ul$ for complex conjugation.
Our first step is to take
two linkages which draw straight lines, as in the
previous section.
In particular, choose $0<a<b$ 
and take $\ul'$ with $\rh{'}A(\conf(\ul'))=[a,b]$,
and $\ul''$ with $\rh{''}B(\conf(\ul''))=[-b,-a]$.
Make sure $\rh{'}A$ and $\rh{''}B$ restricted to
$(a,b)$ and $(-b,-a)$ are analytically trivial covers.
In the \rl\ case we ask that $\rh{'}A$ and $\rh{''}B$ 
be an analytic isomorphisms
 to $[a,b]$ and $[-b,-a]$.

Pick $c>b$ and form
 a new linkage $\ul$ from the disjoint union of $\ul'$ and $\ul''$,
by putting a rigidified square with side length $c$ between $A$ and $B$.

% insert figure 9
\midinsert\vskip 1.8in
\includegraphics{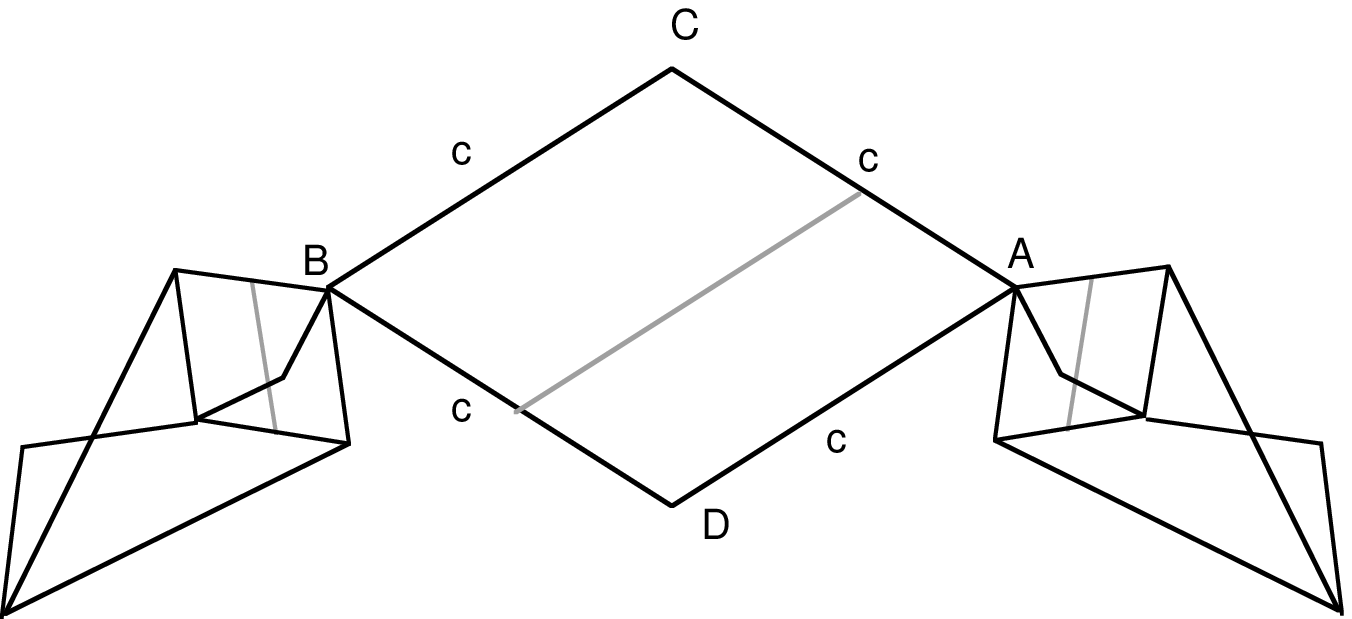}
\capshun {10} { Complex Conjugation}
\endinsert

Note then that if $C$ is the input vertex and $D$ is the output
vertex, then $\ul$ is \qf\ for $z\mapsto \overline z$.
We have an analytic isomorphism 
$$\alpha \colon \conf(\ul')\times \conf(\ul'')\times \{1,-1\}\to \conf(\ul)$$
determined by $\rh{}{\ul'}(\alpha (\varphi ,\psi ,f))=\varphi $,  
$\rh{}{\ul''}(\alpha (\varphi ,\psi ,f))=\psi $, and
$$\alpha (\varphi ,\psi ,f)(C)=(\varphi (A)+\psi (B))/2+f\sqrt{-1}\sqrt{c^2-(\varphi (A)-\psi (B))^2/4}$$
$$\alpha (\varphi ,\psi ,f)(D)=(\varphi (A)+\psi (B))/2-f\sqrt{-1}\sqrt{c^2-(\varphi (A)-\psi (B))^2/4}$$
Using Lemma 6.3, we see that 
$$\rh{}C(\conf(\ul))=\{z\in {\Bbb C} \mid c\ge |z-a|, c\le |z-b|, c\ge |z+a|, c\le |z+b| \}$$
and moreover $\rh{}C$ restricts to an analytically trivial cover
of $U=\{z\in {\Bbb C} \mid c>|z-a|, c<|z-b|, c>|z+a|, c<|z+b| \}$.
But note that if we choose $a$, $b$, and $c$ large enough
then $U$ contains $|z-z_0|\le r$ for some $z_0$.
For example, we may take $a=4r$, $b=8r$, $c=10r$, and $z_0=8r\sqrt{-1}$.
So $\ul$ is functional for $\overline z$ with \rd\ $|z-z_0|\le r$.

In the \rl\ case, we see that $\ul$ is in fact strongly functional,
since with the choices above,  there is unique $(\varphi ,\psi ,f)$
with $\alpha (\varphi ,\psi ,f)(C)=z_0$, namely $f=1$,  and $\varphi =\rh{'}A^{-1}(6r)$,
and $\psi =\rh{''}B^{-1}(-6r)$.
So $\rh{}C$ restricts to a one fold cover of $|z-z_0|\le r$.

\sec{Miscellaneous Comments}

One could also look at something I call a semiconfiguration space
of a linkage $\ul$.
In a semiconfiguration space, you only record the position of
some of the vertices, but ignore the positions of the rest.
In other words it is a projection of the configuration space
to some coordinate $k$-plane, if you ignore all but $k$ vertices.
Semiconfiguration spaces are semialgebraic sets
  (i.e., finite unions of differences of \qas s).
In another paper \cite{K} I  give a complete characterization of
semiconfiguration spaces.
In particular, any compact semialgebraic set $K\subset {\Bbb C}^m$ 
is the semiconfiguration space
of a linkage.
This generalizes Corollary 3.1.

So for example, take any compact polyhedron $K$ in ${\Bbb C}^m$.
Then $K$ is a semialgebraic set 
(since a simplex is a \qas). 
Consequently, $K$
is the semiconfiguration space
of a linkage. 
Put another way, there is a linkage $\ul$ and $m$
vertices of $\ul$ so that if we look at the image of these
$m$ vertices in ${\Bbb C}^m$, they exactly trace out $K$.

Looking at Lemma 3.3 brings up the question of
how the number of fixed vertices affects the topology of $\conf(\ul)$.
In fact, for any linkage $\ul$, there is a linkage $\ul'$ with only three
fixed vertices so that $\conf(\ul)$ is isomorphic to $\conf(\ul')$.
To see this, let $\ul''$ be the linkage obtained from $\ul$ by
adding three vertices $v_0,v_1,v_2$, fixing them at $0,1$, and $\sqrt{-1}$
respectively, and then for each fixed vertex $v$ of $\ul$,
we add three edges, $vv_0$ of length $|z|$, $vv_1$ of length $|z-1|$
and $vv_2$ of length $|z-\sqrt{-1}|$, where $z\in {\Bbb C}$ is the point where
$v$ is fixed.
(This assumes $z$ is not 0, 1 or $\sqrt{-1}$.
If $z$ is one of these values, we identify $v$ with the corresponding $v_i$.)
Let $\ul'$ be obtained from $\ul''$ by unfixing all vertices except
$v_0,v_1,v_2$.
Then $\ul'\subset \ul''$ and $\ul\subset \ul''$ and by Lemma 3.6, both the maps
$\rh{''}{\ul}$ and $\rh{''}{\ul'}$ are isomorphisms.

Now what about linkages with exactly two fixed vertices?
This was the main focus of \cite{KM}.
Assume first that the images of the two fixed vertices
are different, otherwise we could identify them and
obtain a linkage with just one fixed vertex and isomorphic $\conf(\ul)$.
If $\ul$ has exactly two fixed vertices, fixed at
different points of ${\Bbb C}$, there is an involution $\tau \colon \conf(\ul)\to \conf(\ul)$
given by reflection about the line $T$ through the two fixed points.
Usually this involution is nontrivial, it is only trivial when
all  vertices are forced to lie on the line $T$.
By Lemma 6.4 below,
 this can only occur if $\conf(\ul)$ is a finite number of points.
So we have a restriction on the topology of
$\conf(\ul)$, it is either a point or it supports a nontrivial involution.

If the linkage $\ul$ has exactly one fixed vertex (or if all fixed
vertices are fixed to the same point) then 
$\conf(\ul)$ is isomorphic to $\conf(\ul')\times S^1$ for some
linkage $\ul'$ with two fixed vertices.
The linkage $\ul'$ is obtained from $\ul$ by fixing some other vertex.
The $S^1$ factor comes from rotating a planar realization of $\ul'$
around the image of the fixed vertex of $\ul$.

Finally Lemma 3.3 gives restrictions on the topology of
$\conf(\ul)$ if there are no fixed vertices.

I presume that the conclusion to Lemma 6.4 below could be
sharpened to say that $\conf(\ul)$ is a single point.
I presume also that Lemma 6.4 remains true for \rl s.
But in any case we have:

\proclaim{Lemma 6.4}.
Suppose $\ul$ is a classical linkage so that for some line $T\subset {\Bbb C}$,
we have $\varphi (v)\in T$ for all vertices $v$ of $\ul$ and all $\varphi \in \conf(\ul)$.
Then $\conf(\ul)$ is discrete.

\proof
Suppose not.
%First we will do the proof for classical linkages.
Let $\ul$ be such a linkage with the least number of vertices
so that $\conf(\ul)$ is not discrete.
Then for some vertex $v_0$  of $\ul$, there is a one parameter family
$\varphi _t$ in
$\conf(\ul)$, $t\in (-b,b)$ so that $\varphi _t(v_0)\neq \varphi _0(v_0)$
for all $t\in (0,b)$.
So $\varphi _t(v_0)=\varphi _0(v_0)+\alpha (t)z_0$ for some continuous $\alpha $ and constant
 $z_0$ parallel to $T$.
 
 Let $v_1,v_2,\ldots ,v_k$ be the vertices of $\ul$ so that
there is an edge in $\ul$ between $v_i$ and $v_0$.
Order them so that $\varphi _0(v_i)\le \varphi _0(v_j)$ if $i<j$.
(For convenience we identify $T$ isometricly with ${\Bbb R}$.)
Let $n$ be such that $\varphi _0(v_0)\ge \varphi _0(v_i)$ for all $i\le n$ and
$\varphi _0(v_0)<\varphi _0(v_i)$ for all $i>n$.
Since there are no zero length edges, we must have
 $\varphi _0(v_0)\neq \varphi _0(v_i)$ for all $1\le i\le k$.
We must have $\ell (v_0v_i)=|\varphi _t(v_0)-\varphi _t(v_i)|$.
So by continuity, we must have $\varphi _t(v_i)=\varphi _0(v_i)+\alpha (t)z_0$
 for all $i\le k$ and $t\in (-b,b)$.
 In particular, if $\varphi _0(v_i)=\varphi _0(v_j)$, then $\varphi _t(v_i)=\varphi _t(v_j)$
 for all $t\in (-b,b)$.

Note we must have $k\ge 2$ since if $k=0$, then there are no constraints
on the position of $\varphi (v_0)$ for $\varphi \in \conf(\ul)$,
and if $k=1$, the only constraint is that $\varphi (v_0)$ lie in a certain
circle about $\varphi (v_1)$.
In neither case is $\varphi (v_0)$ forced to be in $T$.

Consider the following linkage $\ul'$.
The vertices of $\ul'$ are all the vertices of $\ul$ except for $v_0$.
However, we identify two vertices $v_i$ and $v_j$ connected
to $v_0$ if $\varphi _0(v_i)=\varphi _0(v_j)$.
There is an edge in $\ul'$ between vertices $v$ and $w$ if either:
\item{1)} There is an edge between $v$ and $w$ in $\ul$.
In this case, $\ell '(vw)=\ell (vw)$, they have the same length.
\item{2)} They are both connected to $v_0$,
i.e., $v=v_i$ and $v=v_j$ for some $1\le i\neq j\le k$.
In this case $\ell '(v_iv_j)$=$|\varphi _0(v_i)-\varphi _0(v_j)|$.

Note that by restricting $\varphi _t$ to $\ul'$ we get a parameterized family
$\varphi '_t\in \conf(\ul')$ so that $\varphi '_t(v_i)\neq \varphi '_0(v_i)$ for $t\in (0,b)$.

Take any $\varphi \in \conf(\ul')$.
We will show that $\varphi (v)\in T$ for all vertices $v$ of $\ul'$.
But considering $\varphi '_t$, we see that $\conf(\ul')$ is not 
discrete, which contradicts minimality of $\ul$.

Suppose that $n\neq 0,k$, and thus $\varphi _0(v_1)<\varphi _0(v_0)<\varphi _0(v_k)$.
For any $1<i\le n$ then:
$$
\eqalign{
\ell '(v_1v_i)+\ell '(v_iv_k)&=(\ell (v_0v_1)-\ell(v_0v_i))+(\ell(v_0v_k)+\ell(v_0v_i))\cr
&=\ell (v_0v_1)+\ell(v_0v_k)=\ell '(v_1v_k)\cr}$$
For any $n<i<k$ then: 
$$
\eqalign{
\ell '(v_1v_i)+\ell '(v_iv_k)&=(\ell (v_0v_1)+\ell(v_0v_i))+(\ell(v_0v_k)-\ell(v_0v_i))\cr
&=\ell (v_0v_1)+\ell(v_0v_k)=\ell '(v_1v_k)\cr}$$
so $\varphi (v_i)$ lies on the line segment between $\varphi (v_1)$ and $\varphi (v_k)$.
Moreover, if we identify this line segment isometricly with
$[\varphi _0(v_1),\varphi _0(v_k)]$ we must have $\varphi (v_i)=\varphi _0(v_i)$.
To be precise, there is a Euclidean motion $\beta $ of ${\Bbb C}$ so that
$\varphi (v_i)=\beta \varphi _0(v_i)$ for all $1\le i\le k$.
So from $\varphi $ we get a $\varphi '\in \conf(\ul)$ given by $\varphi '(v)=\varphi (v)$
for $v\neq v_0$ and $\varphi '(v_0)=\beta \varphi _0(v_0)$.
By our supposition, all $\varphi '(v)$ must lie on the line $T$, so in fact we could take
$\beta $ to be a translation parallel to $T$.
In particular, we see that the vertices of $\ul'$ are forced to lie in $T$.

Now suppose that $n=0$. Then:
$$
\eqalign{
\ell '(v_1v_i)+\ell '(v_iv_k)&=(\ell (v_0v_i)-\ell(v_0v_1))+(\ell(v_0v_k)-\ell(v_0v_i))\cr
&=\ell (v_0v_k)-\ell(v_0v_1)=\ell '(v_1v_k)\cr}$$
On the other hand, if $n=k$, then:
$$
\eqalign{
\ell '(v_1v_i)+\ell '(v_iv_k)&=(\ell (v_0v_1)-\ell(v_0v_i))+(\ell(v_0v_i)-\ell(v_0v_k))\cr
&=\ell (v_0v_1)-\ell(v_0v_k)=\ell '(v_1v_k)\cr}$$
So similarly, we see that the vertices of $\ul'$ are forced to lie in $T$.

\commentout{
The proof for \rl s is similar, but one must be careful in certain places.
In particular 
we no longer have $\varphi _t(v_i)=\varphi _0(v_i)+\alpha (t)z_0$
if the edge $v_0v_i$ is not rigid.
If all edges from $v_0$ are flexible, then at least two cables from $v_0$
have to be stretched to their full length and in opposite directions.
That is, there are $i\le n<j$ so that 
$|\varphi _0(v_k)-\varphi _0(v_0)|=\ell (v_kv_0)$ for $k=i,j$.
Otherwise, we could move $v_0$ off the line $T$.
So $\conf(\ul)$ is the union $\conf(\ul_{ij})$
where $\uk_{ij}$ is obtained from $\ul$ by changing
the edges $v_0v_i$ and $v_0v_j$ to be rigid rather
than flexible.
So at least one of the $\conf(\ul_{ij})$ is not discrete,
and we may replace $\ul$ by it,
and thus may assume that there are two rigid edges from $v_0$.
Likewise, if $v_0$ was in only one rigid edge of $\ul$,
we could make one of the other edges rigid and still get
$\conf(\ul)$ not discrete.
So we may as well assume that at least two edges from $v_0$
are rigid.
% Moreover, if $n\ne 0,k$, $\varphi _0$ takes these
% edges to opposite sides of $\varphi _0 (v_0)$.
}

\commentout{
Now in constructing $\ul'$ we must be a bit more careful also
in constructing the edges $v_iv_j$.
If $v_0v_i$ and $v_0v_j$ are rigid edges, then we let
$v_iv_j$ be a rigid edge as before.
}

\qed

\bigbreak

[AK] S.~Akbulut and H.~King, Topology of real algebraic sets,
MSRI Publ.~25, Springer-Verlag (1992).

[AT] S.~Akbulut and L.~Taylor, A topological resolution theorem,
Publ.~I.H.E.S.~53, (1981), pp.~163-195.

[CR] R.~Courant and H.~Robbins, What is Mathematics?, Oxford Univ.~Press (1941).

[K] H.~King, Semiconfiguration spaces of planar linkages, preprint.

[KM] M.~Kapovich and J.~Millson, Universality Theorems for 
configuration spaces of planar linkages, preprint.

\bye